\documentclass[12pt,letterpaper]{article}
\usepackage[utf8]{inputenc}
\usepackage[english]{babel}
\usepackage{amsmath}
\usepackage{amsfonts}
\usepackage{amssymb}
\usepackage{comment}
\usepackage{graphicx}
\usepackage{verbatim,graphics,amsthm,float,subcaption,multirow}
\usepackage{cite,caption,comment,mathrsfs,mathtools}
\usepackage[left=2cm,right=2cm,top=2cm,bottom=2cm]{geometry}
\newtheorem{example}{Example}[section]
\newtheorem{remark}{Remark}[section]

\usepackage{hyperref}
\usepackage{bm}
\usepackage{authblk}
\usepackage{appendix}
\usepackage{algorithm}
\usepackage{algpseudocode}
\usepackage{tikz}
\usetikzlibrary{matrix,chains,positioning,decorations.pathreplacing,arrows}
\usetikzlibrary{shapes.geometric, arrows}
\tikzstyle{startstop} = [rectangle, rounded corners, minimum width=3cm, minimum height=1cm,text centered, draw=black, fill=red!30]
\tikzstyle{process} = [rectangle, minimum width=1cm, minimum height=1cm, text centered, draw=black, fill=orange!30]
\tikzstyle{process_2} = [rectangle, minimum width=2cm, minimum height=1cm, text centered, text width=4cm, draw=black, fill=orange!30]
\tikzstyle{process_3} = [rectangle, minimum width=2.0cm, minimum height=1cm, text centered, text width=2.2cm, draw=black, fill=orange!30]
\tikzstyle{process_4} = [rectangle, minimum width=1cm, minimum height=1cm, text centered, text width=1.5cm, draw=black, fill=orange!30]
\tikzstyle{arrow} = [thick,->,>=stealth]

\usepackage{marginnote}

\newcommand{\change}[1]{{\color{black} #1}}

\let\oldlambda\lambda

\newcommand{\mtbb}[1]{\mathbb{#1}}

\newcommand{\mbbR}{\mtbb{R}}

\newcommand{\mbbU}{\mtbb{U}}
\newcommand{\mbbV}{\mtbb{V}}

\DeclareMathOperator*{\argmin}{argmin}

\title{
Learning quantities of interest from parametric~PDEs: 
\\
An efficient neural-weighted Minimal Residual approach
}
\author[5]{Ignacio Brevis}
\author[1]{Ignacio Muga}
\author[2,3,4]{David Pardo}
\author[2,3]{Oscar Rodriguez}
\author[5]{Kristoffer G. van der Zee}
\affil[1]{Instituto de Matem\'aticas, Pontificia Universidad Cat\'olica de Valpara\'iso, Chile.}
\affil[2]{Universidad del Pa\'is Vasco, Spain.}
\affil[3]{Basque Center for Applied Mathematics, Spain.}
\affil[4]{Ikerbasque, Spain.}
\affil[5]{School of Mathematical Sciences, University of Nottingham, UK.}
\date{February 14, 2024}
\begin{document}
\maketitle
\noindent\makebox[\textwidth][c]{%
\begin{minipage}[t][2cm][t]{9cm}
    \centering
   This paper is dedicated to Leszek Demkowicz on the occasion of his 70th birthday 
\end{minipage}}
\begin{abstract}
\noindent 
The efficient approximation of parametric PDEs is of tremendous importance in science and engineering.
In this paper, we show how one can train Galerkin discretizations to efficiently learn quantities of interest of solutions to a parametric PDE. The central component in our approach is an efficient neural-network-weighted Minimal-Residual formulation, which, after training, provides Galerkin-based approximations in standard discrete spaces that have accurate quantities of interest, regardless of the coarseness of the discrete space. 
\end{abstract}
\tableofcontents
\section{Introduction}
\label{sec:introduction}
\emph{How to learn a reduced order model for parametric partial differential equations (PDEs) given pairs of parameters and solution observations?}
This question has received significant attention in recent years, where a common theme is the hybrid use of neural networks and classical approximations. The aim of this work is to present a new hybrid methodology, based on a neural-network-weighted Minimal Residual (MinRes) Galerkin discretization, that efficiently learns quantities of interest of the PDE solution.
\par
A currently well-established hybrid paradigm for parametric PDEs is to train a neural network to find the coefficients that define an element in a discrete solution space. Effectively, the neural network provides an approximation to the parameter-to-solution map, i.e., the operator that maps parameters to solutions. This has been well explored in the context of non-intrusive Reduced-Basis methods~\cite{HesUbiJCP2018, FreDedManJSC2021, KutPetRasSchCA2022}, where the discrete space is the span of solution snapshots (reduced basis space); see also~\cite{QiaKraPehWilPD2020, BhaHosKovStuSJCM2021, CohDahDevBOOK-CH2022} for related approaches and discussions. More recent is the idea to simply apply this paradigm to the full discrete finite element space~\cite{GeiPetRasSchKutJSC2021, KhaBalJosSarHegKriGan2021, UriParOmeCMAME2022}, cf.~\cite{LuJinPanZhaKarNMI2021}. A disadvantage in these approaches is that the approximations that are produced are no longer guaranteed to satisfy the Galerkin equations of the underlying weak formulation.%
\footnote{In fact, for non-intrusive reduced basis methods, the key motivation is \change{precisely to \emph{interpolate} instead of \emph{project} the full-order model. The reason is that projection leads to Galerkin equations that may not decouple from the full-order model in case of complex nonlinear problems with non-affine dependence~\cite{HesUbiJCP2018}. } }
\par
Because these approaches provide a complete operator-approximation, significant amounts of data may be required to learn it (i.e., to train the neural network). Much less data is expected to be needed to approximate the parameter-to-\emph{quantity-of-interest} map, which is the operator that maps PDE parameters to output quantities of interest (in $\mathbb{R}$) of the PDE solution.%
\footnote{\change{While the parameter-to-PDE-solution map is an operator from some finite-dimensional subset~$\Lambda\subset \mathbb{R}^\rho$ to some infinite-dimensional function space~$\mathbb{V}$, the parameter-to-quantity-of-interest map is just a real-valued function from~$\Lambda\subset \mathbb{R}^\rho$ to $\mathbb{R}$. To have an accurate neural network approximation for the \emph{operator} is a significant challenge, while for the \emph{function} it is much more manageable.}} This has been explored by, e.g., Khoo, Lu \&~Ying~\cite{KhoLuYinEJAM2021}, who use a neural network to map the PDE parameters directly into such quantities. In that approach, there is however no PDE solution approximation.
\par
The objective of this paper is to presents for parametric PDEs, a deep learning methodology for obtaining Galerkin-based approximations in standard discrete spaces that have accurate quantities of interest, regardless of the coarseness of the discrete space. More precisely, the methodology provides an approximation to the parameter-to-quantity-of-interest map, by delivering Galerkin-based approximations on a fixed discrete space for a neural-network-weighted weak formulation. During the (supervised) training, a discrete method is learned by optimizing a weight-function within the formulation. This guarantees that its Galerkin approximations are tailored towards the desired quantities of interest. 
\par
\begin{figure}[h!]
\begin{center}
  \begin{subfigure}[b]{0.32\textwidth}
    \includegraphics[width=\textwidth, height=121px]{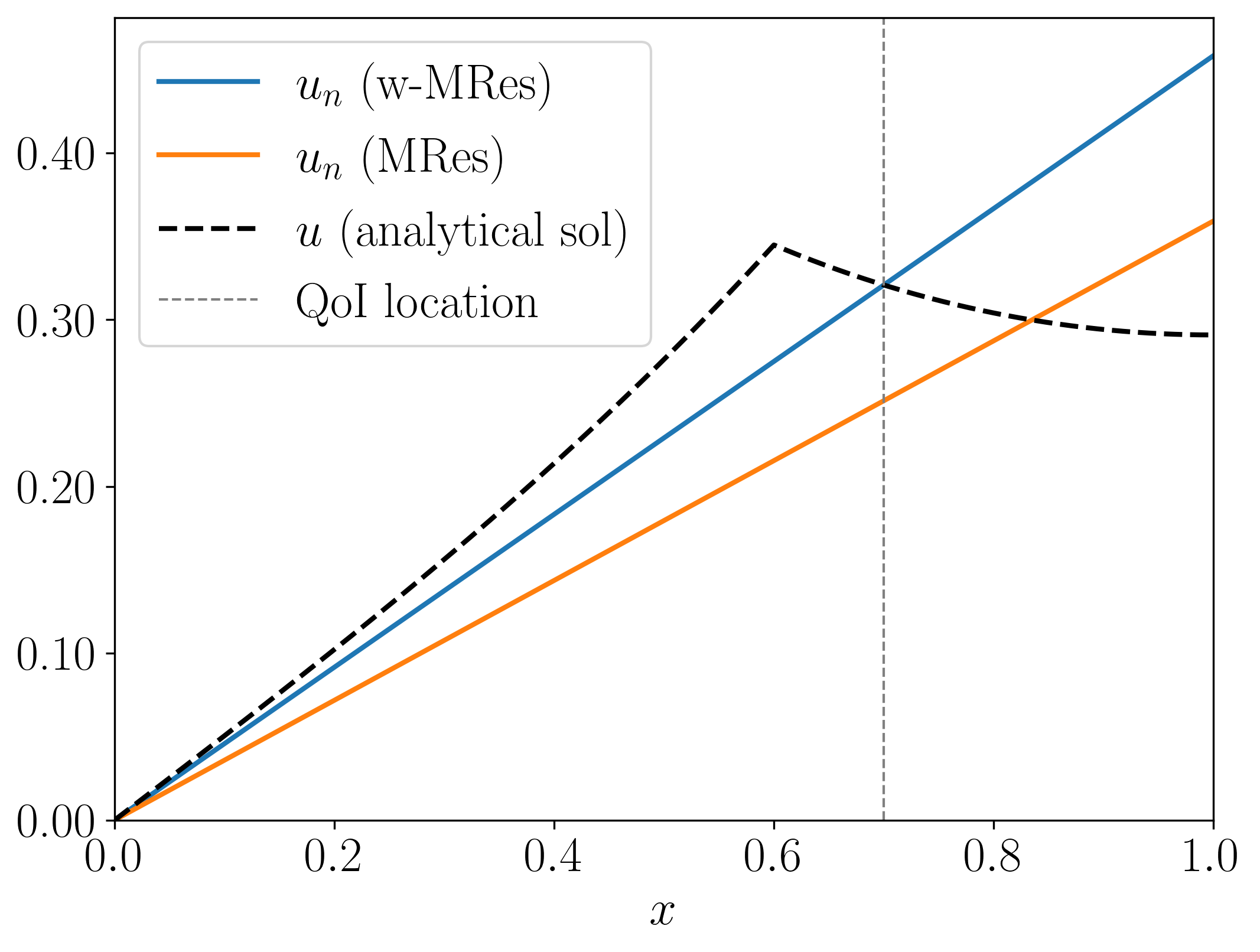}
    \caption{results for $\oldlambda=1.5$}
    \label{sfig:ex01_a}
  \end{subfigure}
  \hfill
  \begin{subfigure}[b]{0.32\textwidth}
    \includegraphics[width=\textwidth, height=121px]{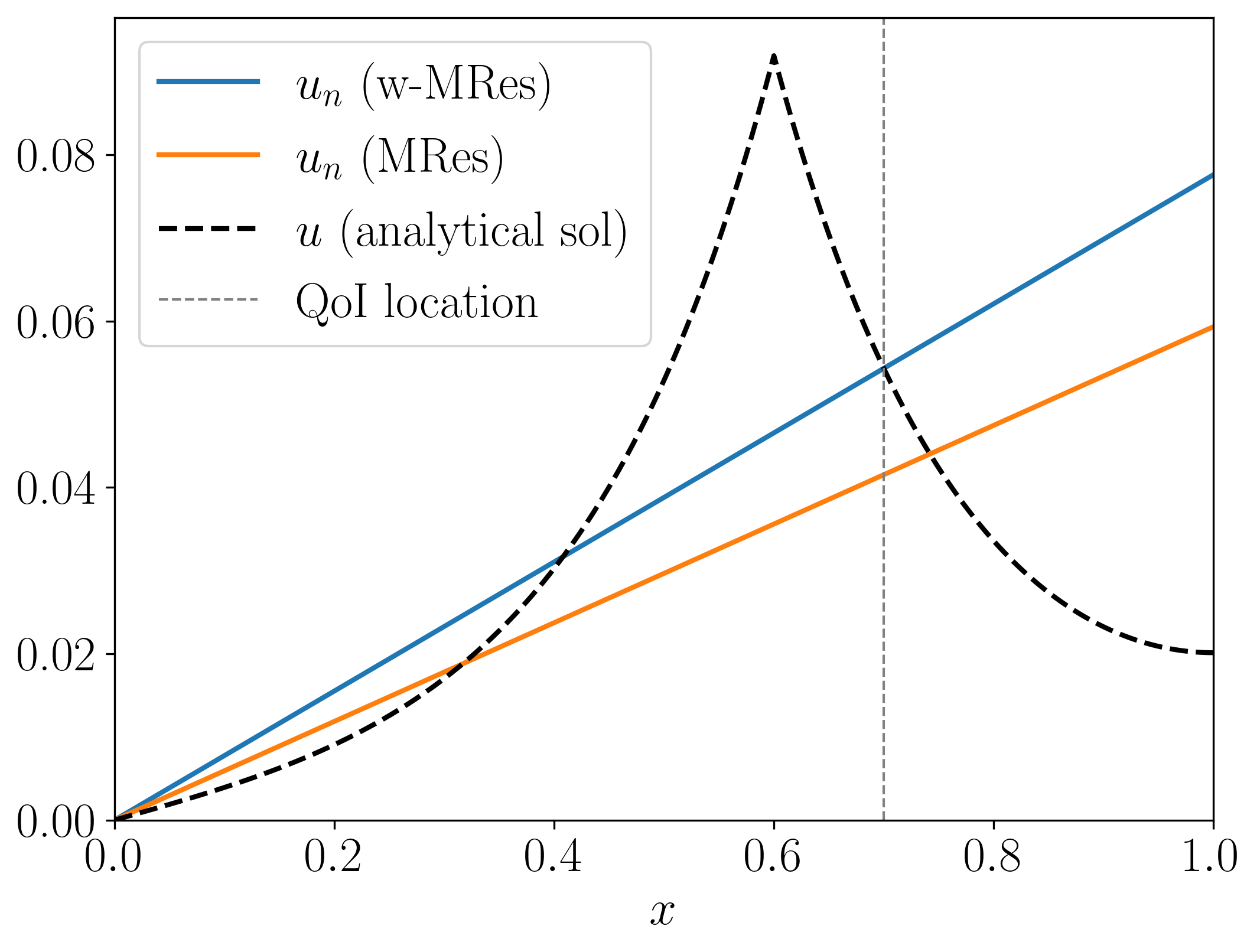}
    \caption{results for $\oldlambda=5.5$}
    \label{sfig:ex01_b}
  \end{subfigure}
  \hfill
  \begin{subfigure}[b]{0.32\textwidth}
    \includegraphics[width=\textwidth, height=121px]{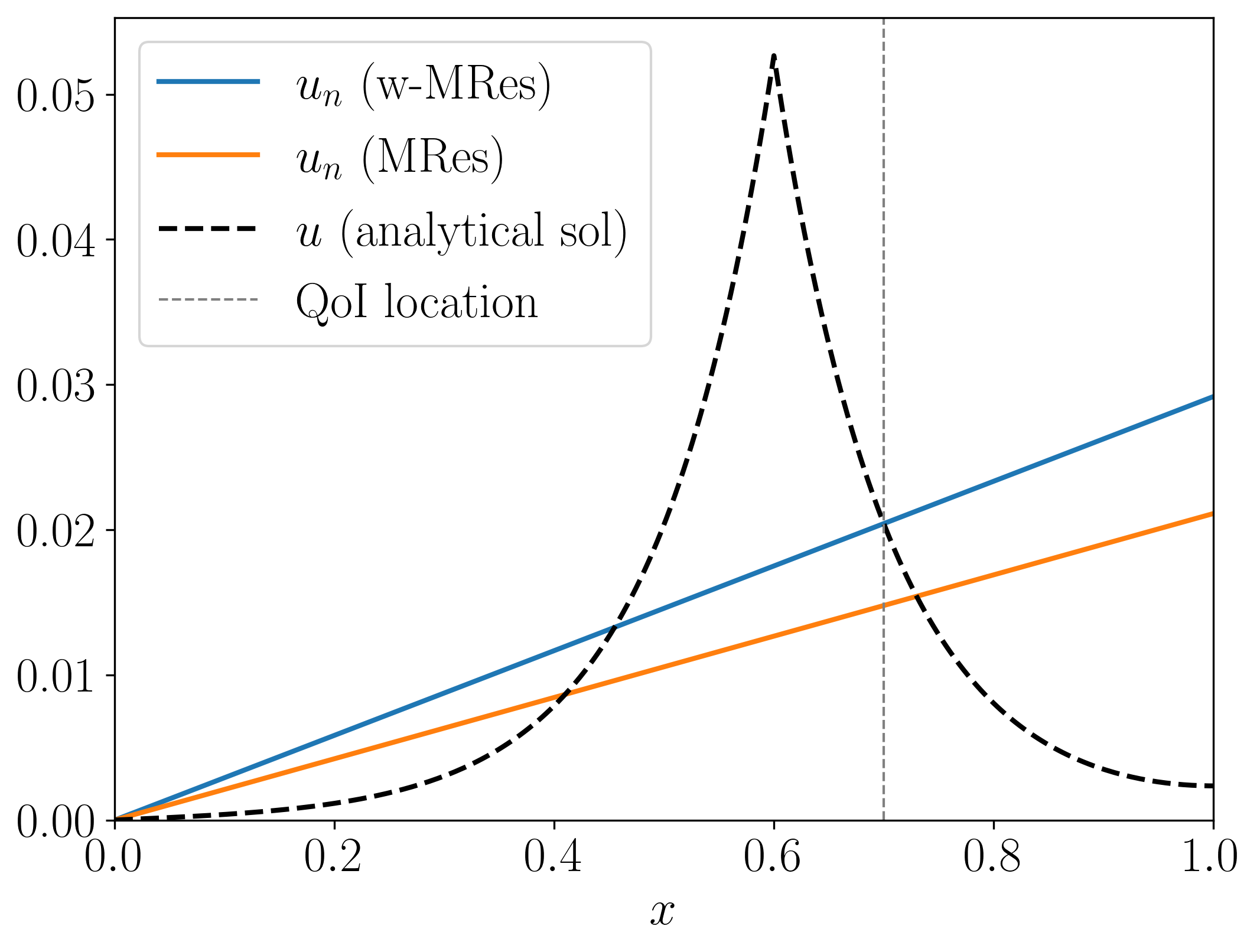}
    \caption{results for $\oldlambda=9.5$}
    \label{sfig:ex01_c}
  \end{subfigure}
\end{center} 
  \caption{\change{Comparison between a standard (unweighted) Galerkin method (MRes) and our trained neural-weighted-method (w-MRes) on a coarse mesh consisting of one linear finite element. The analytical solution is also shown. Note the accuracy of w-MRes at the quantity of interest (QoI), i.e., the value of~$u$ at $x_0=0.7$.}} 
  \label{fig:introductory_example}
\end{figure}

\change{Let us illustrate what our methodology can achieve with a simple example. Consider the following differential equation parametrized by $\oldlambda\in\mathbb R$:
\begin{equation}
\left\{
\begin{array}{rll}
-{d^2\over dx^2}u + \oldlambda^2 u= & \delta_{x_0} & \hbox{in } (0,1), \\
u(0) = & 0, & \\
{du\over dx}(1)= & 0, &
\end{array}
\right.
\label{eq:1d_diff-reac_example}
\end{equation}
where we have a point source at $x_0 = 0.6$, but we are only interested in being accurate at $x=0.7$, even for a coarse mesh. In Figure~\ref{fig:introductory_example}, we approximate the solution of~\eqref{eq:1d_diff-reac_example} with a linear function in only one element (i.e., only one degree of freedom). The standard (unweighted) Galerkin method (indicated by MRes in Figure~\ref{fig:introductory_example}) produces an approximation that is far from being accurate at the quantity of interest location ($x=0.7$); while the neural-network-weighted method (w-MRes) has a $\oldlambda$-dependent trained weight able to recover a discrete solution that is highly accurate at the quantity of interest, for different values of $\oldlambda$.
}
\par
\change{To ensure that the employed weighted  method is stable and consistent, we employ a weighted \emph{Minimal Residual} (abbreviated by MinRes or MRes) weak formulation. The weight function appears in the inner product in the method, and is defined by a neural network. We note that the} theory of Minimal Residual formulations has its origin in the discontinuous Petrov--Galerkin method developed by Demkowicz and Gopalakrishnan~\cite{DemkowiczGopalakrishnan2011, DemGopECM2017}.
\par
Our methodology can be viewed as one where a numerical method is optimized amongst a family of methods, and which utilizes a neural network to learn the optimal method. These ideas go back to Mishra~\cite{MisMINE2018}, who proposed such learning for families of finite differences methods for~PDEs. Similar work on optimizing numerical methods can be found in, e.g.~\cite{RayHesJCP2018, TasZinDedMINE2023}. 
Our current work is a continuation of our previous works~\cite{BreMugZeeCAMWA2021} and \cite{NeuralControl2022}, where we presented a weighted MinRes formulation for the case of parameter dependence in the right-hand-side only, and where we analyzed the well-posedness and convergence in the non-parametric case, respectively.
\par 
The main contributions in our current work are as follows. 
Firstly, we propose the neural-weighted MinRes approach for an abstract and general parametric PDE.
Secondly, we propose new weight functions that depend on the neural network, which lead to piecewise-weighted inner products and allow the practical assembly needed in the training stage.\footnote{And suitable for automatic differentiation in, e.g., Tensorflow.} 
We demonstrate fully-efficient off-line and on-line procedures in the case of affine parameter dependence. 
Thirdly, we propose an adaptive generation of the training set to allow for greedy selection of data. This ensure efficient convergence of errors in the quantity of interest. 
\change{In summary, these contributions provide an efficient methodology for learning the weight function in a weak formulation so that its Galerkin approximations have accurate quantities of interest even on very coarse meshes.}
\par
The contents of the paper are as follows:
Section~\ref{sec:methodology} presents the deep learning methodology for the neural-weighted MinRes formulation.
Section~\ref{sec:efficiency} presents details that demonstrate the efficiency of all involved matrix assemblies, the off-line procedure and on-line procedure. 
Section~\ref{sec:adaptive} presents the adaptive training set approach. 
Section~\ref{sec:num_results} contains numerical experiments. Finally, Section~\ref{sec:concl} discusses conclusions.
\section{Deep learning MinRes methodology}
\label{sec:methodology}

\subsection{The abstract parametric problem}

Let $\mbbU$ and $\mbbV$ be infinite-dimensional Hilbert spaces, with their respective dual spaces $\mbbU^{*}$ and $\mbbV^{*}$. \change{For any $f\in\mathbb V^*$ and $v\in\mbbV$, we will adopt the duality paring notation $\left<f,v\right>_{\mbbV^*,\mbbV}$ to denote the evaluation of $f$ at $v$.} Given a bounded set of parameters\change{\footnote{\change{The dimension of the set of parameters will be denoted by $\rho\in\mathbb N$. In our numerical experiments, we have considered $\rho=1,2$. Higher values of $\rho$ affect the amount of training data needed and therefore the RAM constraints one would have.}}} $\Lambda\subset \mbbR^{\rho}$, we consider parametric bounded invertible linear operators $B(\lambda): \mbbU \to \mbbV^{*}$, and continuous linear functionals $\ell(\lambda)\in\mbbV^*$, where $\lambda \in \Lambda$.  Consequently, we define the parametric family of problems: 
\begin{equation}
\label{eq:abstract_pde}
\left\{
\begin{array}{l}
\hbox{Given $\lambda \in \Lambda$, find $u(\lambda)\in \mbbU$ such that}\\ 
B(\lambda) u(\lambda) = \ell(\lambda) \text{ in }\mbbV^{*}.
\end{array}\right.
\end{equation}
Problems like~\eqref{eq:abstract_pde} are commonly found in variational formulations of parametric PDEs, in which case the linear operator $B(\lambda)$ is defined through a parametrized bilinear form $b(\,\cdot\,,\,\cdot\,;\lambda):\mbbU\times \mbbV \to \mbbR$ such that $\langle B(\lambda)u, v \rangle_{\mbbV^{*},\mbbV} := b(u,v;\lambda)$, for all $(u,v)\in \mbbU\times\mbbV$.  

In many applications, the interest is not focused on the whole solution $u(\lambda)$ of problem~\eqref{eq:abstract_pde}, but in a quantity of interest $q(u(\lambda))$, where $q:\mbbU\to\mbbR$ is a given functional of interest\footnote{For instance, the evaluation of the solution at a single point as in~\eqref{eq:1d_diff-reac_example}, or the average of the solution in some portion of its domain.} (linear or non-linear). In such situations, we are more interested in reproducing the mapping 
$$
\Lambda\ni\lambda \mapsto  q(u(\lambda))\in\mbbR, \quad\hbox{where } u(\lambda) \text{ solves~\eqref{eq:abstract_pde}}.
$$

\subsection{Main idea of Neural-weighted MinRes}

\change{
\subsubsection{Weighted inner-products}
We assume that the Hilbert test space $\mbbV$ comes with an inner-product $(\cdot,\cdot)_\mbbV$.
We want to modify such inner-product by introducing $\lambda$-dependent weights $\omega(\lambda)$, without changing the topology of $\mbbV$\footnote{\change{This is particularly needed to preserve the well-posed structure of problem~\eqref{eq:abstract_pde}}}. The weighted inner-products, denoted by $(\cdot,\cdot)_{\omega(\lambda)}:\mbbV\times\mbbV\to\mbbR$, must satisfy:
\begin{equation}\label{eq:inner-product}
(v,v)_\mbbV \, C_1 \le (v,v)_{\omega(\lambda)} \le (v,v)_\mbbV \, C_2\,,
\end{equation}
for some constants $C_1>0$ and $C_2>0$, uniformly on $v\in\mbbV$ and $\lambda\in\Lambda$.
\begin{example}\label{ex:L2inner}
In an $H_0^1$ context (i.e., $\mbbV=H_0^1(\Omega)$ for $\Omega\subset\mbbR^d$), we can use a positive weight function $\omega(\lambda)\in L^\infty(\Omega)$ and define the weighted $H_0^1$ inner-product
\begin{equation}\label{eq:L2inner}
(f,g)_{\omega(\lambda)} := \int_\Omega \omega(\lambda)\,\nabla f\cdot\nabla g\,.
\end{equation}
It is easy to see that if there are positive constants such that $C_1\le \omega(\lambda)\le C_2$, then~\eqref{eq:inner-product} is satisfied.
\end{example}
The idea will be to train a particular weight function to produce desirable behaviors in our discrete solutions.
}

\change{
\subsubsection{Weighted residual minimization}
We adopt a weighted residual minimization approach to approximate $u(\lambda)\in\mbbU$ solution of~\eqref{eq:abstract_pde}. That is, given a conforming discrete subspace $\mathbb U^n\subset\mbbU$ of dimension $n$, we aim to find 
$u^n(\lambda)\in\mathbb U^n$ such that
\begin{equation}\label{eq:fullresmin}
u^n(\lambda) := \argmin_{w^n\in\mbbU^n}\|\ell(\lambda)-B(\lambda)w^n\|_{\mathbb V^*(\lambda)}\,,
\end{equation}
where the dual norm is now defined using the weighted inner-product, i.e., 
$$
\|\cdot\|_{\mbbV^*(\lambda)}:=\sup_{v\in\mbbV}{\left<\,\cdot\,,v\right>_{\mbbV^*,\mbbV}\over\sqrt{(v,v)_{\omega(\lambda)}}}. 
$$
It is well-known (see,~e.g.,~\cite{CohDahWel_M2AN2012}) that~\eqref{eq:fullresmin} is equivalent to the following mixed (saddle-point) formulation:
\begin{equation}
\left\{ \begin{array}{l@{}ll}
\text{Find } r(\lambda)\in\mbbV \hbox{ and } u^n(\lambda)\in \mbbU^n \text{ s.t.} && \\\\
\big(r(\lambda),v\big)_{\omega(\lambda)} + \big<B(\lambda)u^n(\lambda),v\big>_{\mbbV^*,\mbbV} & = \big<\ell(\lambda),v\big>_{\mbbV^*,\mbbV}\,, & \quad\forall\, v \in \mbbV,\\\\
\big<B(\lambda)w^n,r(\lambda)\big>_{\mbbV^*,\mbbV} & = 0 , & \quad\forall\, w^n \in \mbbU^n.
\end{array} \right.
\label{eq:mixed_formulation}
\end{equation}
Notice that $r(\lambda)\in\mbbV$ is the Riesz representative of the residual $[\ell(\lambda)-B(\lambda)u^n(\lambda)]\in \mbbV^*$ using the weighted inner-product $(\cdot,\cdot)_{\omega(\lambda)}$.

Formulation~\eqref{eq:mixed_formulation} still is infinite-dimensional in the test space $\mbbV$. To come up with a practical scheme, we need to consider a conforming discrete subspace $\mathbb V^m\subset\mbbV$ of dimension $m>n$\footnote{\change{Moreover, to ensure well-posedness and quasi-optimality of the fully-discrete problem, the spaces $\mbbU^n$ and $\mbbV^m$ must satisfy a {\it Fortin} compatibility; see, e.g.,~\cite[Section~4.2]{MugZee_SINUM2020}.}} Thus, the practical (fully-discrete) counterpart of~\eqref{eq:mixed_formulation} will be
\begin{equation}
\left\{ \begin{array}{l@{}ll}
\text{Find } r^m(\lambda)\in\mbbV^m \hbox{ and } u^n(\lambda) \in \mbbU^n \text{ s.t.} && \\\\
\big(r^m(\lambda),v^m\big)_{\omega(\lambda)} + \big<B(\lambda)u^n(\lambda),v^m\big>_{\mbbV^*,\mbbV} & = \big<\ell(\lambda),v^m\big>_{\mbbV^*,\mbbV}\,, & \quad\forall\, v^m \in \mbbV^m,\\\\
\big<B(\lambda)w^n,r^m(\lambda)\big>_{\mbbV^*,\mbbV} & = 0 , & \quad\forall\, w^n \in \mbbU^n.
\end{array} \right.
\label{eq:discrete_mixed_formulation}
\end{equation}}
\begin{remark}
\change{
    Notice that $u^n(\lambda)$ solving~\eqref{eq:discrete_mixed_formulation} is not necessarily the same as $u^n(\lambda)$ solving~\eqref{eq:mixed_formulation} or~\eqref{eq:fullresmin} (they just share the same notation to avoid excessive sub- or supra-scripts). Indeed, $u^n(\lambda)$ of~\eqref{eq:discrete_mixed_formulation} solves the following residual minimization problem in the discrete-dual norm:
    \begin{equation}
u^{n}(\lambda) = 
\argmin_{w^n\in \mbbU^n} \sup_{v^m\in\mbbV^m}
\frac{\big<\ell(\lambda) - B(\lambda) w^n,v^m\big>_{\mbbV^*,\mbbV}}
{\sqrt{(v^m,v^m)_{\omega(\lambda)}}}\,,
\label{eq:discrete_min_res}
    \end{equation}
while $r^m(\lambda)\in \mbbV^m$ is the Riesz representative of the residual $[\ell(\lambda)-B(\lambda)u^n(\lambda)]$ (restricted to $\mbbV^m$) using the weighted inner-product $(\cdot,\cdot)_{\omega(\lambda)}$.}
We refer the reader to~\cite{MugZee_SINUM2020} for a concise and detailed explanation of the equivalences between~\eqref{eq:discrete_mixed_formulation} and~\eqref{eq:discrete_min_res}, and also with the Petrov--Galerkin method with optimal test functions.
\end{remark}

The main idea of neural-weighted MinRes is to work with a trainable weight function $\omega(\lambda)$, that can be trained to minimize errors of $u^n(\lambda)$ in the quantity of interest, i.e., 
\begin{equation}\label{eq:q_min}
{1\over2}\left|q\big(u^n(\lambda)\big)-q\big(u(\lambda)\big)\right|^2\to \min.
\end{equation}
The details about the training process can be found  in Section~\ref{sec:training}.

\subsection{Training of neural-weighted MinRes discrete system}\label{sec:training}

Due to their well-known expressivity properties~\cite{guhring2020expressivity,pmlr-v70-raghu17a}, we choose artificial neural networks~\cite{Goodfellow-DL-2016} to train the weight function $\lambda\mapsto\omega(\lambda)$. 
\change{Therefore, we write $\omega(\lambda;\theta)$ to make explicit the dependence of the weight function in some trainable parameters denoted hereafter by $\theta\in\Theta$.} 
More details on how we choose the weight function to ensure efficiency can be found in the next Section~\ref{sec:efficiency}.

Consider now sets of basis function such that $\mbbU^n = \operatorname{span}\{\varphi_1,\dots,\varphi_n\}$ and $\mbbV^m = \operatorname{span}\{\psi_1,\dots,\psi_m\}$. Thus, system~\eqref{eq:discrete_mixed_formulation} \change{can} be rewritten as:
\begin{equation}
\left\{ \begin{array}{l}
\text{Find }  \change{\bm{r}^{m}(\lambda):=(r_{1}^{m}(\lambda),r_{2}^{m}(\lambda),\dots,r_{m}^{m}(\lambda)) \in \mbbR^{m}}  \\
\change{\text{and } \bm{u}^{n}(\lambda):=(u_{1}^{n}(\lambda),u_{2}^{n}(\lambda),\dots,u_{n}^{n}(\lambda)) \in \mbbR^{n} \text{ s.t.}} \\\\
\arraycolsep=1.4pt\def\arraystretch{1.2}
\left[\begin{array}{c|c}
\mathbb G(\lambda;\theta) & \mathbb B(\lambda) \\
\hline
\mathbb B(\lambda)^{\top} & 0 
\end{array}
\right]
\left[\begin{array}{c}
\change{\bm{r}^{m}(\lambda)} \\ \hline \change{\bm{u}^{n}(\lambda)} 
\end{array}
\right]
= 
\left[\begin{array}{c}
\change{\bm{\ell}(\lambda)} \\ \hline 0
\end{array}
\right], 
\end{array} \right.
\label{eq:discrete_mixed_formulation_system}
\end{equation}
where \change{the solutions $r^m(\lambda)$ and $u^n(\lambda)$ of~\eqref{eq:discrete_mixed_formulation} are recovered by performing the expansions:}
$$r^m(\lambda)  = \sum_{i=1}^{m} \change{r_{i}^{m}(\lambda)} \psi_i \quad\text{and} \quad u^n(\lambda) = \sum_{j=1}^{n} \change{u_{j}^{n}(\lambda)} \varphi_j.$$
Notice that 
\begin{itemize}
    \item $\mathbb G(\lambda;\theta)$ is the Gram matrix obtained from the evaluation of the weighted inner product in the test-space basis functions, i.e., $(\mathbb G(\lambda;\theta))_{ik} = (\psi_k,\psi_i)_{\omega(\lambda;\theta)}$;
    \item $\mathbb B(\lambda)$ denotes the $m\times n$ matrix obtained from the evaluation of the operator $B(\lambda):\mbbU\to\mbbV^*$ in the basis functions, i.e.,  $(\mathbb B(\lambda))_{ij} = \langle B(\lambda) \varphi_j,\psi_i\rangle_{\mbbV^*, \mbbV}$; and
    \item $\bm{\ell}(\lambda)$ represents the load vector computed from the evaluation of the test-space basis functions into the right-hand functional $\ell(\lambda)\in\mbbV^*$, i.e., $(\bm{\ell}(\lambda))_i = \langle \ell(\lambda),\psi_i\rangle_{\mbbV^*, \mbbV}$.
\end{itemize}
\change{Additionally, whenever the quantity of interest functional $q:\mbbU\to\mbbR$ is linear}, $q(u^n(\lambda))$ can be computed in the following way:
\begin{equation}\label{eq:linear_QoI}
    q(u^n(\lambda)) = q\left( \sum_{j=1}^{n} u^n_j(\lambda) \varphi_j\right) =  \sum_{j=1}^{n} u_j^n(\lambda) q\left(\varphi_j\right) = \bm{u}^n(\lambda) \cdot \bm{q},
\end{equation}
where $\bm{q} = (q(\varphi_1),\dots, q(\varphi_n))$ \change{is a vector that can be pre-stored}.

To find an optimal inner-product weight \change{$\omega(\lambda;\theta^*)$}, 
we train \change{an} artificial neural network \change{$\text{NN}(\,\cdot\,;\theta)$} associated with the trainable parameters $\theta\in\Theta$ of the weight. 
\change{First of all, we need to come up with a reliable training set $X_{\text{train}}:=\{\lambda_i,q(u(\lambda_i))\}_{i=1}^{N_s}\subset \Lambda\times\mbbR$, where for each $\lambda_i$
($i=1,\dots,N_s$), the quantity of interest $q(u(\lambda_i))$ (also referred to as \emph{label}; see, e.g. \cite{Goodfellow-DL-2016, Murphy-ML-2012}) corresponds to the evaluation of $q$ at the exact solution of problem~\eqref{eq:abstract_pde} with $\lambda=\lambda_i$.
The training set can be generated in several ways depending on the nature of the physical problem. For example, from measurements or observations, or it can be computed using high-precision numerical schemes.}
The training procedure is executed by minimizing the following cost functional:
\begin{equation}
\mathcal{L}\left(\theta;\change{X_{\text{train}}}\right) := \frac{1}{N_s}\sum_{i=1}^{N_s}\frac{1}{2}\left|\frac{ q(u^n(\lambda_i)) - q(u(\lambda_i))}{q(u(\lambda_i))}\right|^2,
\label{eq:loss_function}
\end{equation}
\change{where the dependence on $\theta$ of the right-hand side of~\eqref{eq:loss_function} is implicit in each $u^n(\lambda_i)$ by means of equation~\eqref{eq:discrete_mixed_formulation_system} with $\lambda=\lambda_i$}.

\begin{figure}[!h]
\begin{center}
\begin{tikzpicture}[node distance=1.5cm]
%
\node (mix_form) [startstop] {\begin{minipage}[t]{6.7cm} \vspace*{-0.3cm}$$ 
\begin{array}{l}
\text{Find } (\bm{r}^m(\lambda), \bm{u}^n(\lambda)) \in \mbbR^{m}\times \mbbR^{n}, \text{s.t.}\\
\arraycolsep=1.4pt\def\arraystretch{1.2}
\left[\begin{array}{c|c}
\mathbb{G}(\lambda;\theta) & \mathbb{B}(\lambda) \\
\hline
\mathbb{B}(\lambda)^{\top} & 0 
\end{array}
\right]
\left[\begin{array}{c}
\bm{r}^m(\lambda) \\ \hline \bm{u}^n(\lambda) 
\end{array}
\right]
= 
\left[\begin{array}{c}
\bm{\ell}(\lambda) \\ \hline 0
\end{array}
\right]
\end{array}
$$ \vspace*{-0.3cm}\end{minipage}};
\node (omega) [process_2, above of=mix_form, yshift=1.0cm] {inner-product weight \\ $\omega(\lambda;\theta) 
$};
\node (lambda) [process, left of=mix_form, xshift=-3.0cm] {$\lambda$
};
\node (train_set) [process_3, left of=lambda, xshift=-0.8cm] {$\mathrm{training\; set}$ $\left\{\lambda,q(u(\lambda))\right\}$ 
};
\node (qoi) [process_4, right of=mix_form, xshift=3.4cm] {$q({u}^n(\lambda))$
};
\node (theta) [process, right of=omega, xshift=3.1cm] {$\text{NN}(\lambda;\theta)$};
\node (loss) [process, right of=qoi, xshift=1.2cm] {loss function};
\node (qoi_ts) [process_4, below of=mix_form, yshift=-0.6cm] {$q(u(\lambda))$
};
\draw [arrow] (theta) -- (omega);
\draw [arrow] (omega) -- (mix_form);
\draw [arrow] (lambda) -- (mix_form);
\draw [arrow] (mix_form) -- (qoi);
\draw [arrow] (qoi) -- (loss);
\draw [arrow] (train_set)-- (lambda);
\draw [arrow] (train_set) |- (qoi_ts);
\draw [arrow] (qoi_ts) -| (loss);
\draw [arrow] (loss) |- (theta);
\draw [dashed] (-8.2,-2.9) rectangle (9.0,1.6);
\draw (-7.3,-2.6) node {Epoch};
\draw (6.5,2.2) node {update $\theta$};
\end{tikzpicture}
\caption{Supervised training process for the artificial neural network and the inner-product weight.}
\label{fig:Training_scheme}
\end{center}
\end{figure}
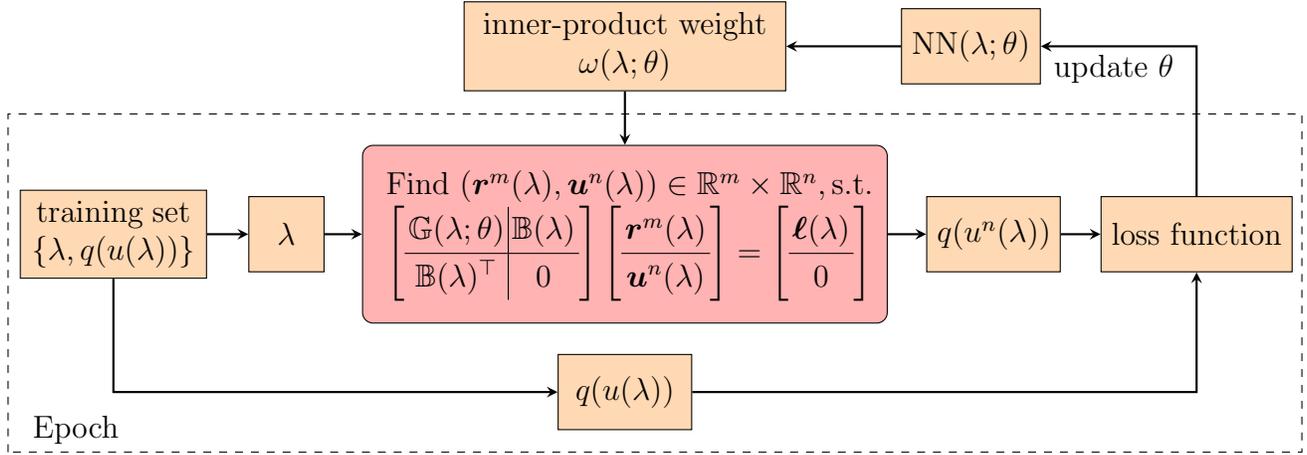

Figure \ref{fig:Training_scheme} summarizes the training procedure, describing the process at each epoch. That is, given $\theta\in\Theta$, we build the weight function $\omega(\lambda;\theta)$. For each $\lambda_i$ in the training set, we solve the system~\eqref{eq:discrete_mixed_formulation_system} to compute $u^n(\lambda_i)$ using the weighted inner product associated with the weight $\omega(\lambda_i;\theta)$. Finally, the cost functional~\eqref{eq:loss_function} is used to update the value of $\theta$, and repeat. 

\subsection{Exemplification}
\change{
To put the above abstract formulations into context, let us consider the introductory example~\eqref{eq:1d_diff-reac_example}, with $\oldlambda\in\Lambda:=[1,10]\subset \mathbb{R}$. 
Here, we use the following trial and test spaces:
$$
\mbbU=\mbbV=H_{(0}^1(0,1):=\left\{v\in L^2(0,1): {dv\over dx}\in L^2(0,1) \hbox{ and } v(0)=0\right\}.
$$
Thus, the linear quantity of interest is given by $q(u)= u(0.7)$, and the involved operators are:   
$$\langle B(\oldlambda)u,v\rangle_{\mbbV^*,\mbbV} := \int_{0}^{1} {du\over dx}{dv\over dx} + \oldlambda^2
\int_{0}^{1} u v  \quad\text{ and }\quad \langle\ell(\oldlambda),v\rangle_{\mbbV^*,\mbbV}:= v(x_0).$$ 
We use a positive weight function $\omega(\oldlambda)\in L^\infty(\Omega)$ and define a weighted $H_{(0}^1(0,1)$ inner-product as follows:
\begin{equation}\label{eq:H1inner}
(v,\nu)_{\omega(\oldlambda)} := \int_{0}^{1} \omega(\oldlambda) \left({dv\over dx}{d\nu\over dx} + v\nu\right)\,.
\end{equation}

Given a uniform partition of the interval $(0,1)$ into $m$ elements, we can consider conforming piecewise linear finite-element spaces such that $\mbbU^n\subsetneq\mbbV^m$. A simple setting for a trainable $\omega(\lambda;\theta)$ is to use piecewise constant functions on each of those elements. Therefore, an inner product for the above system can be expressed as:
\begin{equation}\label{eq:pwc_inner}
(v,\nu)_{\omega(\lambda;\theta)} := \sum_{i=1}^{m}c_i(\oldlambda;\theta)\int_{(i-1)/m}^{i/m}  \left({dv\over dx}{d\nu\over dx} + v\nu\right),
\end{equation}
for some trainable scalars $\oldlambda\mapsto \{c_1(\oldlambda;\theta),...,c_m(\oldlambda;\theta)\}$.
In the next Section, we seek to employ this type of inner product as they can be easily pre-assembled.}
\section{Ensuring Efficiency}
\label{sec:efficiency}

\subsection{Efficient assembly: Precompute sub-domain Gram matrices}
\label{sec:Gram}

We will work with variational formulations of PDEs. Hence, 
our inner products must depend on the PDE domain variable, say $x\in\Omega\subset\mbbR^d$. \change{Here we consider $\mbbV$ as the Hilbert space $H^{s}(\Omega)$, with $s\in\mathbb N$,  endowed with the inner product
$$(v,\nu)_{\mbbV} := \sum_{|\alpha|=0}^{s}\left(D^{\alpha} v,D^{\alpha} \nu\right)_{L^2(\Omega)},$$
where $\alpha=(\alpha_1,...,\alpha_d)$ is a multi-index and 
$$D^{\alpha} v := \frac{\partial^{|\alpha|} v}{\partial_{x_1}^{\alpha_1}\cdots\partial_{x_n}^{\alpha_d}}.$$} 
Moreover, we want to use conforming piecewise polynomial approximation spaces \change{$\mbbU^n\subset\mbbU$ and $\mbbV^m\subset\mbbV$. The trial space may be coarse for approximating the PDE solution, but it has to be fine enough for its functions to reach the desired values at the quantities of interest, for any $\lambda\in\Lambda$.}  

To assemble the Gram matrix $\mathbb G(\lambda;\theta)$ in~\eqref{eq:discrete_mixed_formulation_system}, we need to compute the inner products $(\mathbb G(\lambda;\theta))_{ik} = (\psi_k,\psi_i)_{\omega(\lambda;\theta)}$ for every value of $\lambda$ in the training set, for every trainable parameter $\theta$ under consideration, and for every pair of piecewise polynomial test-space basis functions in $\{\psi_1,...,\psi_m\}$. This could be an extremely expensive task for inner products like~\eqref{eq:L2inner} or~\eqref{eq:H1inner}. To overcome these difficulties, we propose the following strategies: 
\begin{itemize}
\item Let $\Omega_h$ be a partition of $\Omega$, defined as a collection of finitely many open connected elements $T\subset \Omega$ with Lipschitz boundaries, such that  $\overline{\Omega}=\cup_{T\in\Omega_h}\overline{T}$.
\item Define a discrete test space $\mathbb V^m\subset\mathbb V$ as a conforming piecewise polynomial space over the mesh $\Omega_h$.
\item Define a (coarser) discrete trial space $\mathbb U^n\subset\mbbU$ as a conforming piecewise polynomial space over the mesh $\Omega_h$.
\item Define a non-overlapping domain decomposition of $\Omega$, consisting of $n_a$ patches $\{\Omega_l\}_{l=1}^{n_a}$ of elements of $\Omega_h$, such that $\cup_{l=1}^{n_a}\overline{\Omega_l}=\overline{\Omega}$.
\item \change{Define a trainable parametric family of weighted inner products in $\mbbV$ as follows:
\begin{equation}\label{eq:weighted_inner_prod}
(v,\nu)_{\omega(\lambda;\theta)} := \sum_{|\alpha|=0}^{s}\left(\omega(\lambda;\theta) D^{\alpha} v,D^{\alpha} \nu\right)_{L^2(\Omega)}= \sum_{|\alpha|=0}^{s} \sum_{l=1}^{n_a} c_{l}(\lambda;\theta)\left(D^{\alpha} v,D^{\alpha} \nu\right)_{L^2(\Omega_l)},
\end{equation}
where $\omega(\lambda;\theta):= c_{1}(\lambda;\theta)\chi_{\Omega_1} + \cdots + c_{n_a}(\lambda;\theta)\chi_{\Omega_{n_a}}$ with positive coefficients $c_l$, and $\chi_{\Omega_l}$ being the characteristic function of $\Omega_l$.} 
\item In terms of the PDE parameter $\lambda\in\Lambda$, \change{the positive constants $c_l(\lambda;\theta)$ will be obtained from a neural network $\text{NN}:\Lambda \to \mbbR^{n_a}$, whose input parameter is $\lambda$, and whose output parameters are the $n_a$ coefficients 
$\{c_{l}(\lambda;\theta)\}_{l=1}^{n_a}$}. 
\end{itemize}
With this setting, the computation of coefficients $(\mathbb G(\lambda;\theta))_{ik}$ of the Gram matrix becomes a piecewise polynomial integration, where we can use quadrature rules to perform exact computations (avoiding the uncertain problem of integrating the neural network; see, e.g.,~\cite{rivera2021quadrature}). However, the main advantage is that we can pre-assemble the Gram matrix before training in order to avoid integration routines during the training procedure. 
\begin{remark}\label{rem:alternative}
\change{As an alternative, we can use neural networks whose inputs are $(\lambda,l)$, and whose output is only the coefficient $c_{l}(\lambda;\theta)$ for the patch $\Omega_l$. This may be challenging to implement when using some machine learning packages, as there would be no label for each training set element $(\lambda,l)$ in this setting. In fact, each label $q(u(\lambda))$ would be associated to a set of elements $\{(\lambda,1),\dots,(\lambda,n_a)\}$ in the training set.}
\end{remark}

Consider the following explanatory example.
\begin{example}\label{ex:G_assemble}
In the same context of Example~\ref{ex:L2inner}, let $\mathbb V^m:=\operatorname{span}\{\psi_1,...,\psi_m\}\subset H_0^1(\Omega)$. In this case, using the  expression~\eqref{eq:weighted_inner_prod}, the coefficients of the Gram matrix become
\begin{alignat*}{2}
(\mathbb G(\lambda;\theta))_{ik}=(\psi_k,\psi_i)_{\omega(\lambda;\theta)}= & \sum_{l=1}^{n_a}
\int_{\Omega_l}
\omega(\lambda;\theta)\,\nabla\psi_k\cdot\nabla\psi_i\\
= & \sum_{l=1}^{n_a}\change{c_{l}(\lambda;\theta)}
\underbrace{\int_{\Omega_l}\,\nabla\psi_k\cdot\nabla\psi_i}_{=:(\mathbb M^{l})_{ik}},
\end{alignat*}
where $(\mathbb M^{l})_{ik}$ corresponds to one coefficient of a (sparse) sub-domain Gram matrix $\mathbb M^{l}$, which can be pre-assembled before training using exact quadrature rules. Therefore, at each Epoch and for any $\lambda$, the full Gram matrix can be rapidly assembled using the formula
$$
\mathbb G(\lambda;\theta)  =\sum_{l=1}^{n_a}\change{c_{l}(\lambda;\theta)}
\,\mathbb M^{l}.
$$
\end{example}

\begin{remark}[Practical \change{positive} piecewise polynomial weights]
\change{In terms of the PDE variable $x\in\Omega$, one may generally use piecewise polynomial weights 
$\omega(\lambda;\theta)\in L^\infty(\Omega)$ such that $\omega(\lambda;\theta)\big|_{\Omega_l}$ is a positive polynomial of higher order (instead of a constant). However, when using nodal basis functions on simplices, describing the set of all positive polynomials of degree $p>1$ is on another level of complexity (not so for polynomials of degree $p=1$).} 
\end{remark}

\subsection{Efficient assembly: Precompute stiffness matrices and load vectors with unit parameters}\label{sec:Stiff}
Recall from Section~\ref{sec:training} that we are also working with a $m\times n$ matrix
$\mathbb B(\lambda)$ resulting from the evaluation of the operator $B(\lambda):\mbbU\to\mbbV^*$ in the basis functions of the discrete spaces $\mathbb U^n$ and $\mathbb V^m$. That is, $(\mathbb B(\lambda))_{ij} := \langle B(\lambda) \varphi_j,\psi_i\rangle_{\mbbV^*, \mbbV}$
where $\operatorname{span}\{\varphi_1,\dots,\varphi_n\}=\mbbU^n$ and $\operatorname{span}\{\psi_1,\dots,\psi_m\}=\mbbV^m$.

We further assume that the operator $B(\lambda)$ allows an affine decomposition of the form 
$$
B(\lambda) = B_0 + \sum_{l=1}^{n_b}\Phi_l(\lambda)\, B_l\,,
$$
where $\{B_l:\mathbb U\to \mathbb V^*\}_{l=0}^{n_b}$ are continuous linear operators (independent from $\lambda$) and $\{\Phi_l:\Lambda\to \mathbb R\}_{l=1}^{n_b}$ are given functions of the set of parameters.
Naturally, this will induce a self-explained matrix decomposition of the form 
\begin{equation}
    \mathbb B(\lambda)= \mathbb B_0 +\sum_{l=1}^{n_b}\Phi_l(\lambda)\,\mathbb B_l\, ,
    \label{eq:B_affine}
\end{equation}
where the matrices $\{\mathbb B_l\}_{l=0}^{n_b}$ can be pre-assembled before training following the same idea of the previous Section~\ref{sec:Gram}.

The same assumption is made on the parametrized right-hand side functional 
$\ell(\lambda)\in \mathbb V^*$. Namely, express 
$\ell(\lambda) = \ell_0 + \sum_{l=1}^{n_l} \Psi_l(\lambda)\, \ell_l$; next define loading vectors
\begin{equation}
\bm{\ell}(\lambda) = \bm{\ell}_0 + \sum_{l=1}^{n_l} \Psi_l(\lambda)\, \bm{\ell}_l\,;
\label{eq:L_affine}
\end{equation}
and proceed as before pre-assembling $\lambda$-independent vectors before training.

\change{
\subsection{Offline and Online procedures}}

\change{For the offline procedure,} once we have performed all the pre-assemblies described in previous Sections~\ref{sec:Gram} and~\ref{sec:Stiff}, we start defining the neural network and its training settings, i.e., the neural network architecture (which contains the trainable parameters) and the training algorithm with its \change{non-trainable} hyperparameters. 
We observe two main possible architectures (cf.~Eq.\eqref{eq:weighted_inner_prod} and Remark~\ref{rem:alternative}):
\begin{itemize}
\item
$\lambda \to \text{NN}(\lambda;\theta) \to \{c_{l}(\lambda;\theta)\}_{l=1}^{n_a}$; and
\item 
$(\lambda,l) \to \text{NN}(\lambda,l;\theta) \to 
c_{l}(\lambda;\theta)$.
\end{itemize}

Defining a training set $\change{X_{\text{train}}}=\{\lambda_i,q(u(\lambda_i))\}_{i=1}^{N_s}$ (and a validation set if necessary), we train the neural network by minimizing the loss function~\eqref{eq:loss_function}. This loss will deliver the trained parameter of the neural network, i.e.,
\change{
\begin{equation}\label{eq:trained}
\theta^* \approx \displaystyle\argmin_{\theta\in \Theta} \mathcal{L}\left(\theta,X_{\text{train}}\right).
\end{equation}}
We emphasize that the cost functional $\mathcal{L}$ is constrained by the linear system~\eqref{eq:discrete_mixed_formulation_system}.  
Additionally, in Section~\ref{sec:adaptive} we describe an adaptive training procedure that can be used to improve the performance of the trained neural network.

\change{In the online procedure, we take the} trained parameter $\theta^*$, \change{and} for any given $\lambda\in\Lambda$
we solve the system: 
\begin{equation}
\arraycolsep=1.4pt\def\arraystretch{1.2}
\left[\begin{array}{c|c}
\mathbb G(\lambda;\theta^*) & \mathbb B(\lambda) \\
\hline
\mathbb B(\lambda)^{\top} & 0 
\end{array}
\right]
\left[\begin{array}{c}
\bm{r}^n(\lambda) \\ \hline \bm{u}^n(\lambda) 
\end{array}
\right]
= 
\left[\begin{array}{c}
\bm{\ell}(\lambda) \\ \hline 0
\end{array}
\right], 
\end{equation}
\change{which implies solving a linear system of size $n+m$. Finally, we compute $q(u^n(\lambda))$\footnote{\change{Alternatively, we can use the formula $\bm{u}^n(\lambda)\cdot\bm{q}$ for linear quantities of interest; see~\eqref{eq:linear_QoI}.}} as an approximation of the exact quantity of interest $q(u(\lambda))$.}
\section{Neural-Network Adaptive Training Set}
\label{sec:adaptive}

In classical training for machine learning models, we usually rely on a fixed amount of observations, which are split between training and validation sets (sometimes also on a test set). However, with the increasing number of machine learning models, different approaches to training were developed. Those techniques range from avoiding overfitting from large amounts of data to avoiding underfitting by generating new data from small training sets. 
The training time and final result will depend, among other elements, on the election of an adequate training set. 
Some well-known methodologies of adaptive training set are bootstrapping in neural networks \cite{franke2000bootstrapping}, important sampling \cite{nabian2021efficient}, and mesh refinement in PINNs \cite{wu2023comprehensive}. We implement a methodology to add samples to training and validation sets at specific training epochs.

\subsection{Adaptive generation of training data}
\label{subsec:adaptive_training}

We employ an adaptive training step to improve the training process without significantly increasing computational cost. Based on the idea of adaptive integration presented in \cite{rivera2021quadrature}, we define an adaptive training set that will iteratively incorporate the elements of the validation set with the worst performance (according to certain criteria). 
Particularly in this work, we include in the training set those elements of the validation with a loss evaluation greater than $\gamma$ times the loss of the training set. Algorithm~\ref{alg:cap} describes this idea. 

\begin{figure}[!h]
\centering
\includegraphics[scale=1]{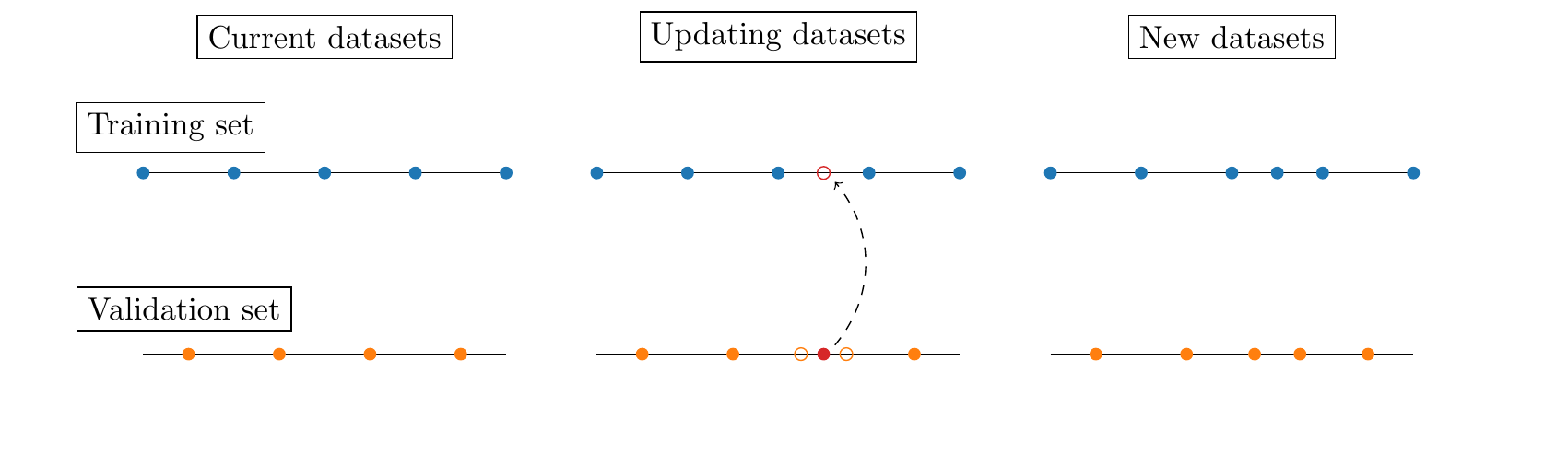}
\caption{Adaptive training and validation sets.}
\label{fig:adaptive_datasets}
\end{figure}

\begin{algorithm}
\caption{Adaptivity for the training set}\label{alg:cap}
\begin{algorithmic}
\Require $X_{\text{train}}$, $X_{\text{val}}$, $\gamma$ \Comment{\change{$X_{\text{train}} = \{\lambda_i,q(u(\lambda_i))\}_{i=1}^{N_s}$, $X_{\text{val}} = \{\lambda_j,q(u(\lambda_j))\}_{j=1}^{N_v}$}}
\While{stage $<$ $\#$stages} \Comment{$\#$stages: number of stages to update the training set}
    \While{epoch $<$ \#epochs} \Comment{$\#$epochs: number of epochs}
        \State $\theta^* \gets \argmin_{\theta} \mathcal{L}(\theta; X_{\text{train}})$ \Comment{neural network training}
    \EndWhile
    \State $\mathcal{L}_{\text{train}} \gets \mathcal{L}(\theta^*; X_{\text{train}})$
    \For{$x_{\text{val}} \in X_{\text{val}}$}
        \If{$\mathcal{L}(\theta^*; \{x_{\text{val}}\}) > \gamma \mathcal{L}_{\text{train}}$} \Comment{Chosen criteria}
            \State $X_{\text{train}} \gets X_{\text{train}} \cup \{x_{\text{val}}\}$
        \EndIf
    \EndFor
    \State update $X_{\text{val}}$ \Comment{$X_{\text{val}}$ is updated when new points are added to $X_{\text{train}}$}
\EndWhile
\end{algorithmic}
\end{algorithm}
\begin{remark}
If $\gamma$ is greater than one, we avoid overfitting by adding the validation-set points with a loss value that lags behind the training-set loss value.
\end{remark}
\change{The above adaptive algorithm keeps track of the neural network performance during the training process by computing the relative squared error (given by the loss function~\eqref{eq:loss_function}) of each element in the validation set and comparing them with the loss function of the training set.} 
We propose as the validation set the middle points (center of mass in higher dimensions) of the elements in the training \change{ set from an equispaced grid}. So each time that a validation-set point is included as part the training set, we update the validation set adding the new middle point from the training set. Figure~\ref{fig:adaptive_datasets} illustrates the training and validation set update procedure.
For example, 
\begin{enumerate}
    \item Select a training set from the coarse grid of the parameter space $\Lambda$.
    \item We define the validation set as the middle points of the coarse grid. 
    \item Run a usual minimization method (SDG, Adam, RMSProp, etc) to train the artificial neural network with the above training set for a certain number of epochs.
    \item Evaluate the loss function within the elements of the validation set.
    \item Add the validation-set points with the worst results to the training set (according to certain criteria, see Algorithm~\ref{alg:cap}).
    \item Update the validation set by removing the points added to the training set and adding the new needed middle points \change{of the training set}.
\end{enumerate}

\section{Numerical Results}
\label{sec:num_results}

We present some numerical examples to show the different features of the method. 
In the next examples, unless otherwise stated, we employ feed-forward neural networks with three hidden layers and ten neurons on each hidden layer. 
The input of the neural network is the PDE parameter $\lambda$, and the input layer of the neural network has $\rho$ neurons (with $\rho := \dim \Lambda$). The output of the neural network is the set of coefficients associated with the weighted inner product. Thus, for a piecewise constant weight function defined on $n_a$ patches, the output layer contains $n_a$ neurons, i.e., 
\begin{equation}
\nonumber
\lambda \longrightarrow \text{NN}(\lambda;\theta) \longrightarrow 
\big\{c_1(\lambda;\theta), \dots, c_{n_a}(\lambda;\theta)\big\}.
\end{equation}
To ensure positive constants, we choose the $\operatorname{softplus}$ function as the activation function for the output layer (while the hyperbolic tangent function is used on every other layer). 
To train the neural network, we employ the Adam optimization algorithm~\cite{DBLP:journals/corr/KingmaB14} within TensorFlow, with decay rates $\beta_1=0.9$ and $\beta_2=0.999$, regularization term for numerical stability $\epsilon=1\text{e}^{-16}$, EMA momentum equal to $0.99$, and learning rate $\alpha = 1.\text{e}^{-4}$. 
Table~\ref{tab:nn-architecture} displays a TensorFlow scheme of the neural network architecture.

\begin{table}[H]
\centering
\begin{tabular}{|l|l|l|}
\hline
\multicolumn{3}{|c|}{\textbf{ML-MinRes neural network architecture}} \\ \hline
\textbf{Layer (type)} & \textbf{Output Shape} & \textbf{Activation Function} \\
\hline \hline
Input non-trainable layer ($\lambda$) & (None, $\rho$) & $\tanh$ \\ \hline
Hidden Dense layer & (None, 10) & $\tanh$ \\ \hline
Hidden Dense layer & (None, 10) & $\tanh$ \\ \hline
Hidden Dense layer & (None, 10) & $\tanh$ \\ \hline
Output Dense layer ($c$) & (None, $n_a$) & $\operatorname{softplus}$ \\ 
\hline
\end{tabular}
\caption{Neural network architecture employed in most of the following examples.}
\label{tab:nn-architecture}
\end{table}

\subsection{1D diffusion-reaction equation with one parameter}
\label{subsec:example_01}

We consider \change{the introductory 1D diffusion-reaction equation~\eqref{eq:1d_diff-reac_example}, where} the right-hand side of this problem is given by a Dirac delta function located at $x_0 = 0.6$,
\change{with $\oldlambda\in\Lambda:= [1,10]$, and} the linear quantity of interest $q(w)= w(0.7)$, for all $w\in\mbbU$.

Notice that the exact solution to problem~(\ref{eq:1d_diff-reac_example}) is:
\begin{equation}
\nonumber
[u(\oldlambda)](x) = \begin{cases} \change{k_1}(\exp(\oldlambda x) - \exp(-\oldlambda x)), & \text{if } x < x_0, \\ \change{k_2}\left(\exp(\oldlambda x) + \exp(\oldlambda(2-x))\right), & \text{if } x > x_0, \end{cases}
\end{equation}
where 
$$
\change{k_1} = \displaystyle \frac{1}{2\oldlambda}\frac{\exp(\oldlambda x_0) + \exp(\oldlambda (2-x_0))}{1 + \exp(2\oldlambda)}\qquad\hbox{ and }\qquad \change{k_2} = \displaystyle\frac{1}{2\oldlambda}\frac{\exp(\oldlambda x_0) - \exp(-\oldlambda x_0)}{1 + \exp(2\oldlambda)}.
$$

Given the continuous trial and test spaces $\mbbU=\mbbV = H_{(0}^{1}(0,1)$, we use the following variational formulation of problem~\eqref{eq:1d_diff-reac_example}:
\begin{equation}
\left\{
\begin{array}{l@{\;}lr}
\multicolumn{3}{l}{\text{Find } u(\oldlambda) \in \mbbU \text{ such that:}} \\
\underbrace{\displaystyle \int_{0}^{1} {d\over dx}u(\oldlambda) {dv\over dx}}_{\left<B_0 u(\oldlambda),v\right>_{\mathbb V^*,\mathbb V}} + \oldlambda^2
\underbrace{\int_{0}^{1} u(\oldlambda) v}_{\left<B_1u(\oldlambda),v\right>_{\mathbb V^*,\mathbb V}} &= v(x_0), & \forall\, v\in \mbbV.
\end{array}
\right.
\end{equation} 
Thus, we define $B(\oldlambda):\mathbb U\to \mathbb V^*$ as $B(\oldlambda):=B_0+\oldlambda^2 B_1$ and $\ell(\oldlambda):=\delta_{x_0}$ (which is independent from $\oldlambda$).

The discrete trial space $\mbbU^n$ will be the space generated by the function $\varphi(x)=x$, while the discrete test space $\mbbV^m$ will consist of conforming piecewise linear functions over a uniform mesh of four elements. Hence, $n=1$ and $m=4$.

\change{The inner product will be the one of expression~\eqref{eq:pwc_inner}. Notice that we can define the weighted inner product in multiple ways. We choose to add one constant per element of the test space mesh. In this way, we can benefit from the richness of the discrete test space to find a good approximation on the coarse space $\mathbb{U}_n$}.

With all these ingredients we construct the mixed formulation:
\change{
\begin{equation}
\label{eq:mixed_exp1}
\nonumber
\left\{ \begin{array}{l@{}ll}
\text{Find } r^{m}(\oldlambda)\in \mbbV^{m} \hbox{ and } u^{n}(\oldlambda)\in\mbbU^{n} \text{ such that for all } v^{m} \in \mbbV^{m} \hbox{ and } w^{n} \in \mbbU^{n} && \\\\
\displaystyle\sum_{i=1}^{4}\!c_i(\oldlambda)\!\!\int_{(i-1)/4}^{i/4}\!\!\left(\! {d \over dx}r^{m}(\oldlambda) {d v^{m}\over dx}+ r^{m}(\oldlambda) v^{m}\!\right)\!+\!\int_{0}^{1}\!\!\frac{d}{dx} u^{n}(\oldlambda)\frac{d v^{m}}{dx}\! +\! \oldlambda^2\!\!
\int_{0}^{1}\!\! u^{n}(\oldlambda) v^{m} &{} = v^{m}(x_0), &\\\\
\displaystyle\int_{0}^{1} \frac{d w^{n}}{dx} \frac{d}{dx} r^{m}(\oldlambda) + \oldlambda^2
\int_{0}^{1} w^{n} r^{m}(\oldlambda) &{} = 0 , & 
\end{array} \right.
\end{equation}
}
and proceed with the methodology described in Section~\ref{sec:methodology}. We use the neural network architecture described in Table~\ref{tab:nn-architecture} and no adaptive training has been performed for this case. The training set is composed of $N_s = 10$ elements, such that $\oldlambda \in \{1, 2, 3, \dots, 10\}$.

\begin{figure}[h!]
\begin{center}
  \begin{subfigure}[b]{0.32\textwidth}
    \includegraphics[width=\textwidth, height=121px]{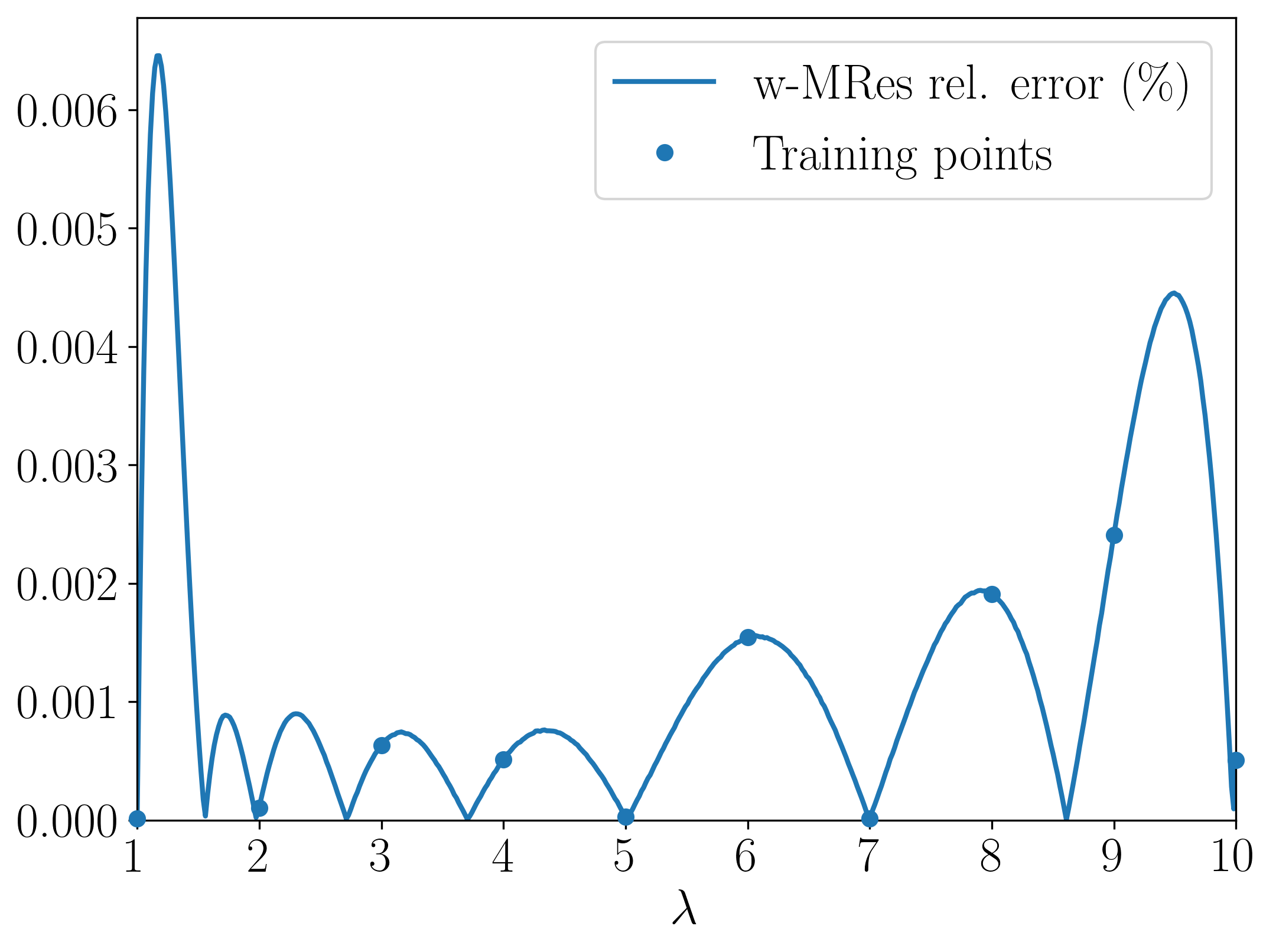}
    \caption{w-MinRes method rel. error}
    \label{sfig:ex01_d}
  \end{subfigure}
  \hfill
  \begin{subfigure}[b]{0.32\textwidth}
    \includegraphics[width=\textwidth, height=121px]{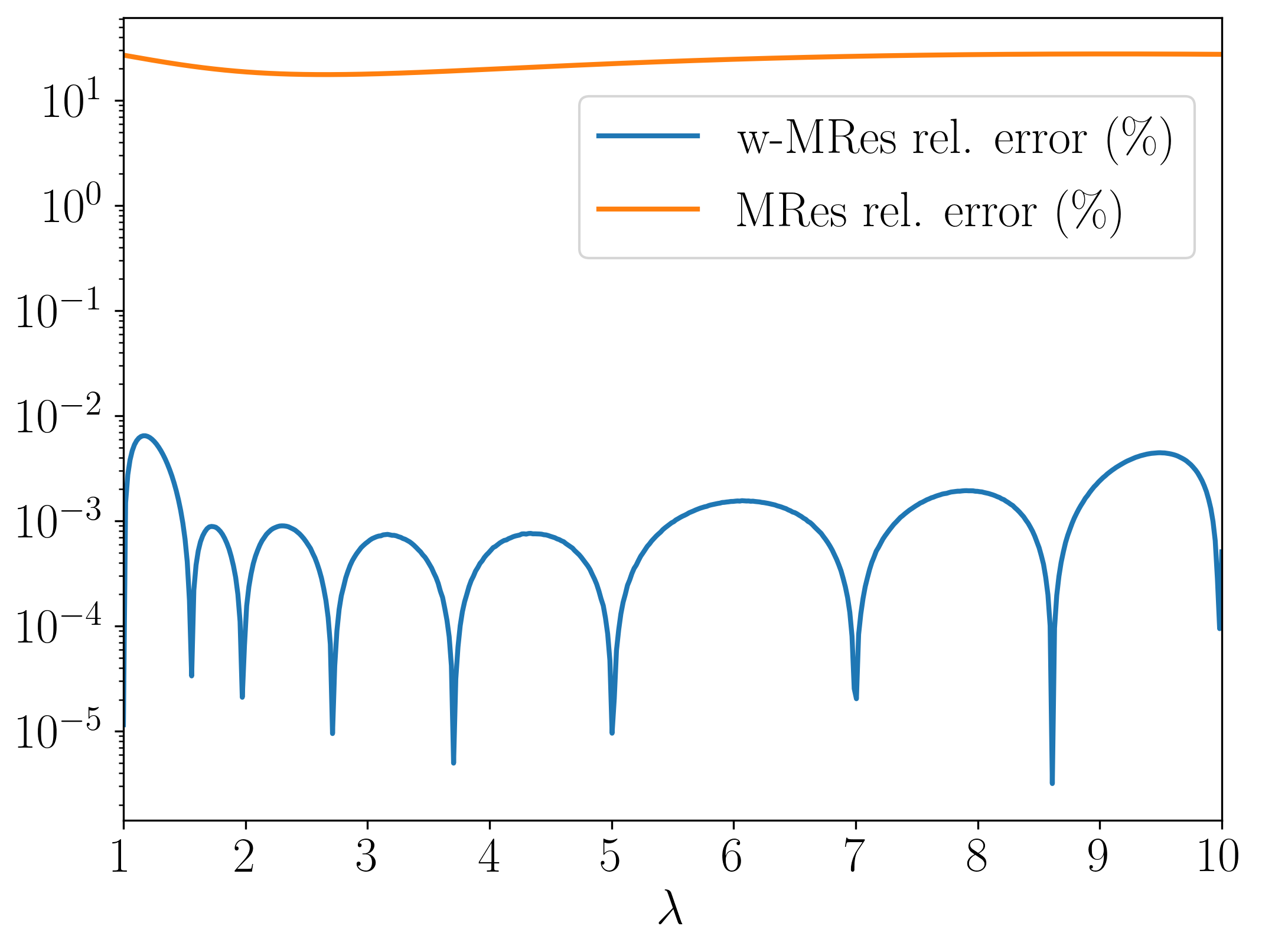}
    \caption{comparison rel. error} 
    \label{sfig:ex01_e}
  \end{subfigure}
  \hfill
  \begin{subfigure}[b]{0.32\textwidth}
    \includegraphics[width=\textwidth, height=121px]{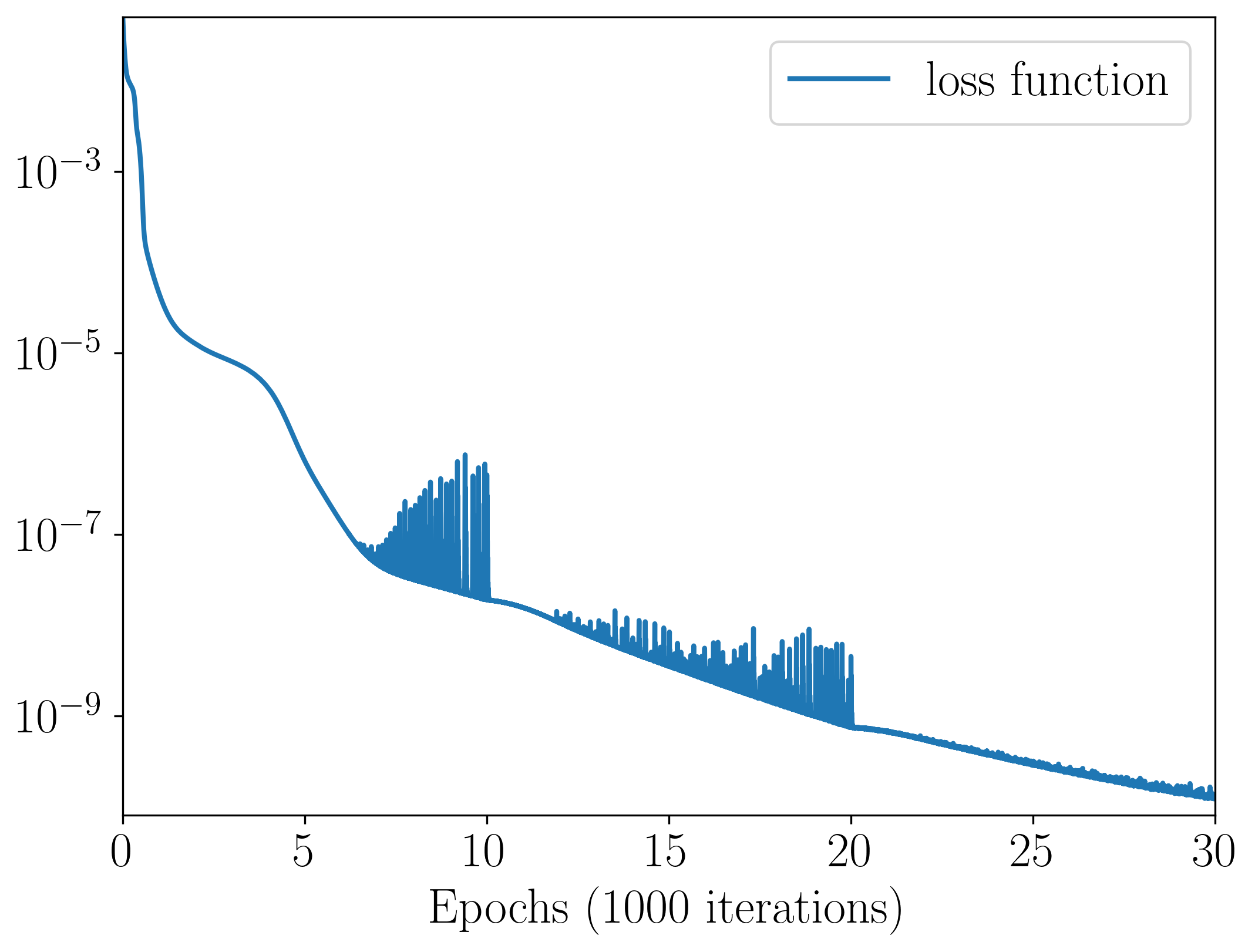}
    \caption{loss function}
    \label{sfig:ex01_f}
  \end{subfigure}
\end{center} 
  \caption{ML-MinRes performance for 1D parametric diffusion reaction example.} 
  \label{fig:example_1_results_some_parameters_trial_one_element}
\end{figure}

\change{Figure~\ref{fig:introductory_example} (in Section \ref{sec:introduction}) and Figure}~\ref{fig:example_1_results_some_parameters_trial_one_element} show the performance of ML-MinRes or weighted-Minres (w-MinRes) method, after training the inner-product weight for the example described above. \change{Figure~\ref{fig:introductory_example} compares} the results for the standard (unweighted) MinRes versus the ML-MinRes methods for $\oldlambda \in \{1.5, 5.5, 9.5\}$. 
Subfigure~\ref{sfig:ex01_e} compares the relative error in logarithmic scale between the MinRes and ML-MinRes methods for a fine-grid test set. 
The last subfigure shows the loss function through epochs. Notice that for this training we employed three different learning rates $[1\text{e}^{-3}, 1\text{e}^{-4}, 1\text{e}^{-5}]$, each one for $10.000$ iterations, in order to avoid oscillations in the loss function. 
\change{The numerical solutions of the w-MinRes method give a good approximation of the quantity of interest (with a relative error of less than $10^{-2}$ \%, see subfigure~\ref{sfig:ex01_e}), even when the discrete trial space is composed of a single element. We can observe this from the intersection of the w-MinRes solution (blue line) with the analytical solution (black dashed line) in the quantity of interest (grey dashed line) in subfigures~\ref{sfig:ex01_a},~\ref{sfig:ex01_b}, and~\ref{sfig:ex01_c}.}

\subsection{1D diffusion-reaction equation with two parameters}

We consider 1D diffusion-reaction problem with a two-dimensional parameter $\lambda = (\alpha,\beta)$. 
The right-hand side of this problem is given by a Dirac delta function located at $x_0 = 0.6$. 
The first component of the parameter $\alpha$ multiplies the diffusion term and the second component $\beta$ multiplies the reaction term. Thus, we have the following problem: 
\begin{equation}
\left\{
\begin{array}{rll}
-\alpha^2 \frac{d^2}{dx^2}u + \beta^2 u & = \delta_{x_0}, & \hbox{in } (0,1). \\
u(0) & = 0, & \\
\frac{du}{dx}(1) & =0, &
\end{array}
\right.
\label{eq:ex3_1d_diff-reac_example}
\end{equation}
with $\lambda = (\alpha,\beta) \in \Lambda:=[1,10]^2$. 
We set the quantity of interest as $q(w)= w(0.7)$, for all $w\in\mbbU$. Notice that the exact solution to this problem~(\ref{eq:ex3_1d_diff-reac_example}) is
\begin{equation}
[u(\lambda)](x) = \begin{cases} \displaystyle \change{k_1}\left(\exp \left(\frac{\beta}{\alpha} x \right) - \exp \left(-\frac{\beta}{\alpha} x \right)\right), & \text{if } x < x_0, \\ \displaystyle \change{k_2}\left(\exp \left(\frac{\beta}{\alpha} x \right) + \exp\left(\frac{\beta}{\alpha}(2-x)\right)\right), & \text{if } x > x_0, \end{cases}
\label{eq:example2_analytic_sol}
\end{equation}
where 
$$\change{k_1} = \frac{1}{2\alpha \beta} \frac{ \displaystyle\exp\left(\frac{\beta}{\alpha} x_0\right) + \exp\left(\frac{\beta}{\alpha} (2-x_0)\right)}{ \displaystyle 1 + \exp\left(2\frac{\beta}{\alpha}\right)}\quad\text{ and }\quad \change{k_2} = \frac{1}{2\alpha\beta}\frac{\displaystyle\exp\left(\frac{\beta}{\alpha} x_0\right) - \exp\left(-\frac{\beta}{\alpha} x_0\right)}{\displaystyle 1 + \exp\left(2\frac{\beta}{\alpha}\right)}.$$

As in the previous example, we set the trial and test spaces to be $\mbbU=\mbbV = H_{(0}^{1}(0,1)$, and the variational formulation is given by
\begin{equation}
\left\{
\begin{array}{l@{\;}lr}
\multicolumn{3}{l}{\text{Find } u(\lambda) \in \mbbU \text{ such that:}} \\
\displaystyle \alpha^2\underbrace{\int_{0}^{1} \frac{d}{dx}u(\lambda) \frac{dv}{dx}}_{\left<B_1u(\lambda),v\right>_{\mbbV^*,\mbbV}} + \,\beta^2\underbrace{\int_{0}^{1} u(\lambda) v}_{\left<B_2u(\lambda),v\right>_{\mbbV^*,\mbbV}} &= \displaystyle v(x_0) , & \forall\, v\in \mbbV.
\end{array}
\right.
\end{equation} 
Thus, we define $B(\lambda):\mathbb U\to \mathbb V^*$ as $B(\lambda):=\alpha^2 B_1+\beta^2 B_2$ and $\ell(\lambda):=\delta_{x_0}$ (which is independent from $\lambda$).
In addition, we use the same inner product together with the same (discrete) trial and test space of the previous example in Section~\ref{subsec:example_01}. \change{The discrete mixed formulation is given by} 

\begin{figure}[h!]
\begin{center}
  \begin{subfigure}[b]{0.32\textwidth}
    \includegraphics[width=\textwidth, height=121px]{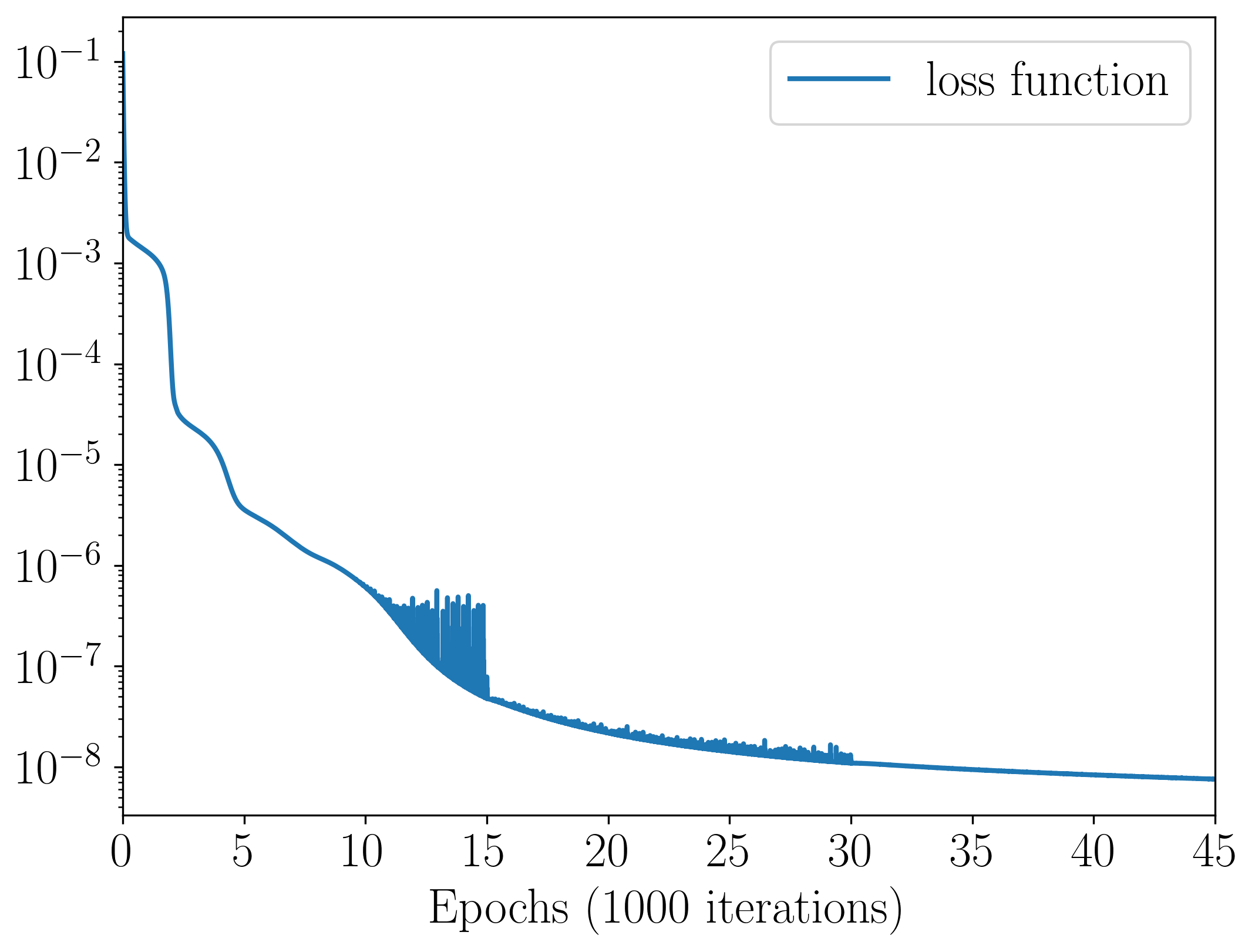}
    \caption{loss function}
    \label{sfig:diff_reac_two_param_loss_func}
  \end{subfigure}
  \begin{subfigure}[b]{0.32\textwidth}
    \includegraphics[width=\textwidth,  height=121px]{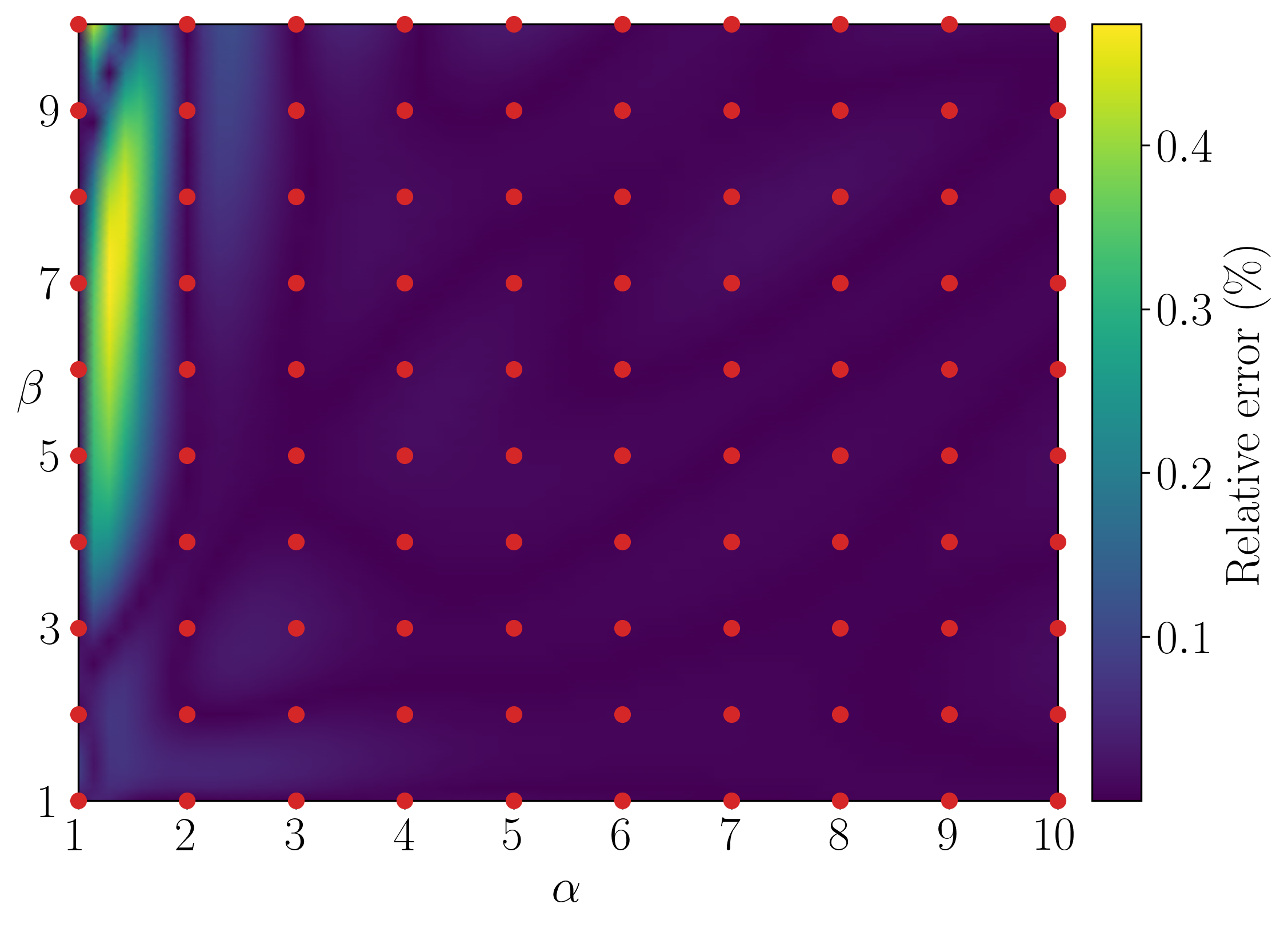}
    \caption{w-MinRes method rel. error}
    \label{sfig:diff_reac_two_param_rel_error}
  \end{subfigure}
\end{center}  
  \caption{ML-MinRes performance for two-parameters diffusion reaction example. Red dots show the training points.}
  \label{fig:example_3_1D_two_params_diff_reaction}
\end{figure}

\change{
\begin{equation}
\nonumber
\left\{ \begin{array}{l@{}l}
\text{Find } r^{m}(\lambda)\in \mbbV^{m} \hbox{ and } u^{n}(\lambda)\in\mbbU^{n} \text{ such that for all } v^{m} \in \mbbV^{m} \hbox{ and } w^{n} \in \mbbU^{n} & \\\\
\displaystyle \sum_{i=1}^{4}c_i(\lambda)\!\int_{(i-1)/4}^{i/4} \hspace{-5px}\left(\frac{d}{dx} r^{m}(\lambda) \frac{d  v^{m}}{dx} + r^{m}(\lambda) v^{m}\right) + \alpha^2\!\int_{0}^{1}\hspace{-1px} \frac{d}{dx} u^{n}(\lambda) \frac{d v^{m}}{dx}  + \beta^2\!
\int_{0}^{1}\hspace{-1px} u^{n}(\lambda) v^{m} &{} = v^{m}(x_0), 
\\\\
\displaystyle \alpha^2\int_{0}^{1} \frac{d w^{n}}{dx}  \frac{d}{dx} r^{m}(\lambda) + \beta^2
\int_{0}^{1} w^{n} r^{m}(\lambda) &{} = 0. 
\end{array} \right.
\end{equation}
}
Notice that in this example, the artificial neural network has two neurons in the input layer (and four in the output layer). For the training set, we select 100 equispaced elements in the set $(1,10)\times(1,10)$ (see red dots in Subfigure~\ref{sfig:diff_reac_two_param_rel_error}). In addition, we employed three different learning rates $[1\text{e}^{-3}, 1\text{e}^{-4}, 1\text{e}^{-5}]$ during the training, each one for $15.000$ iterations. 

Subfigure~\ref{sfig:diff_reac_two_param_loss_func} displays the loss function through epochs for the three different learning rates $[1\text{e}^{-3}, 1\text{e}^{-4}, 1\text{e}^{-5}]$. 
Subfigure~\ref{sfig:diff_reac_two_param_rel_error} shows the relative error in linear scale for a test set (fine-grid in $\Lambda$) and also display the training points in red dots. The picture shows that the relative error is below $0.5\%$ after the training. 

\subsection{1D advection with parametric right-hand side}

Let us consider the parametric 1D advection problem:
\begin{equation}
\left\{
\begin{array}{rll}
\displaystyle\frac{du}{dx} & = \ell({\oldlambda}), & \hbox{in } (0,1), \\
u(0) & = 0, &
\end{array}
\right.
\label{eq:1d_advection_param_rhs}
\end{equation}
where $\left<\ell(\oldlambda),v\right>_{\mbbV^*,\mbbV} :=\int_\oldlambda^1 (x - \oldlambda)v(x)dx$, for all $v\in\mbbV:=L^2(0,1)$ and $\oldlambda \in \Lambda:=[0,1]$. We set the quantity of interest to be $q(w) = w(0.9)$, for all $w\in\mbbU:=H_{(0}^{1}(0,1)$. 

We use the variational formulation,
\begin{equation}
\left\{
\begin{array}{lr}
\text{Find } u(\oldlambda) \in H_{(0}^{1}(0,1) \text{ such that} \\
\underbrace{\displaystyle \int_{0}^{1} \frac{d}{dx} u(\oldlambda) v}_{\left<Bu,v\right>_{\mbbV^*,\mbbV}} = \underbrace{\int_{\oldlambda}^{1}(x-\oldlambda) v(x)dx}_{\left<\ell(\oldlambda),v\right>_{\mbbV^*,\mbbV}}, & \forall\, v \in L^{2}(0,1).
\end{array}
\right.
\label{eq:1d_advection_param_rhs_v-f_s-f}
\end{equation}
Notice that the exact solution to this problem is
\begin{equation}
    \nonumber
    [u(\oldlambda)](x) = \left\{
    \begin{array}{rc}
    (x - \oldlambda)^2, & \textrm{ if } x > \oldlambda, \\
    0, & \textrm{ if } x < \oldlambda.
    \end{array}
    \right.
\end{equation}
In addition, we use \change{a weighted $L^{2}$} inner product \change{with the same weight function of the previous examples}, together with the same one-dimensional discrete trial space $\mbbU^n$. 
However, the discrete test space $\mbbV^m$ will be a piecewise constant space over four uniform elements (same as the weight)\footnote{Notice that in this particular example, there is no need to discretize the test space because we can treat the problem using a weighted least squares approach (see~\cite[Section 3.1]{NeuralControl2022})}. \change{The discrete mixed formulation, in this case, is given by:}
\change{
\begin{equation}
\nonumber
\left\{ \begin{array}{l@{}ll}
\text{Find } r^{m}(\oldlambda)\in \mbbV^{m} \hbox{ and } u^{n}(\oldlambda)  \in \mbbU^{n} \text{ such that} && \\\\
\displaystyle \sum_{i=1}^{4}c_i(\oldlambda)\int_{(i-1)/4}^{i/4} r^{m}(\oldlambda) v^{m} + \int_{0}^{1} \frac{d}{dx} u^{n}(\oldlambda) v^{m} &{} = \displaystyle \int_{\oldlambda}^{1}(x-\oldlambda) v^{m}(x)dx, & \quad\forall\, v^{m} \in \mbbV^{m},\\\\
\displaystyle\int_{0}^{1} \frac{dw^{n}}{dx} r^{m}(\oldlambda) &{} = 0 , & \quad\forall\, w^{n} \in \mbbU^{n}.
\end{array} \right.
\end{equation}
}

For this example, we employ the adaptive training set strategy from Subsection~\ref{subsec:adaptive_training}, with $30000$ epoch per stage and $\gamma=5$.   
We also define the initial training set $X_{\text{train}}=\{(\oldlambda_i,q(u(\oldlambda_i))\}_{i=1}^{N_s}$ composed of $N_s = 11$ elements such that $\oldlambda_{i+1} - \oldlambda_{i} = 0.1$ for all $i \in \{1,2,\dots,N_s-1\}$. Notice that some quantities of interest may be zero, so in this case we modify the loss function, by adding a small regularization term
\begin{equation}
    \nonumber
    \mathcal{L}_{\mathrm{reg}}\left(\theta;X_{\text{train}}\right) := \frac{1}{N_s}\sum_{i=1}^{N_s}\frac{1}{2} \left|\frac{q(u(\oldlambda_i)) - q(u^{n}(\oldlambda_i))}{q(u(\oldlambda_i)) + \epsilon_0}\right|^2, 
\label{eq:loss_function_reg}
\end{equation}
where $\epsilon_0 = 10^{-6}$.
\begin{figure}
    \centering
    \begin{subfigure}[b]{0.32\textwidth}
    \includegraphics[width=\textwidth, height=121px]{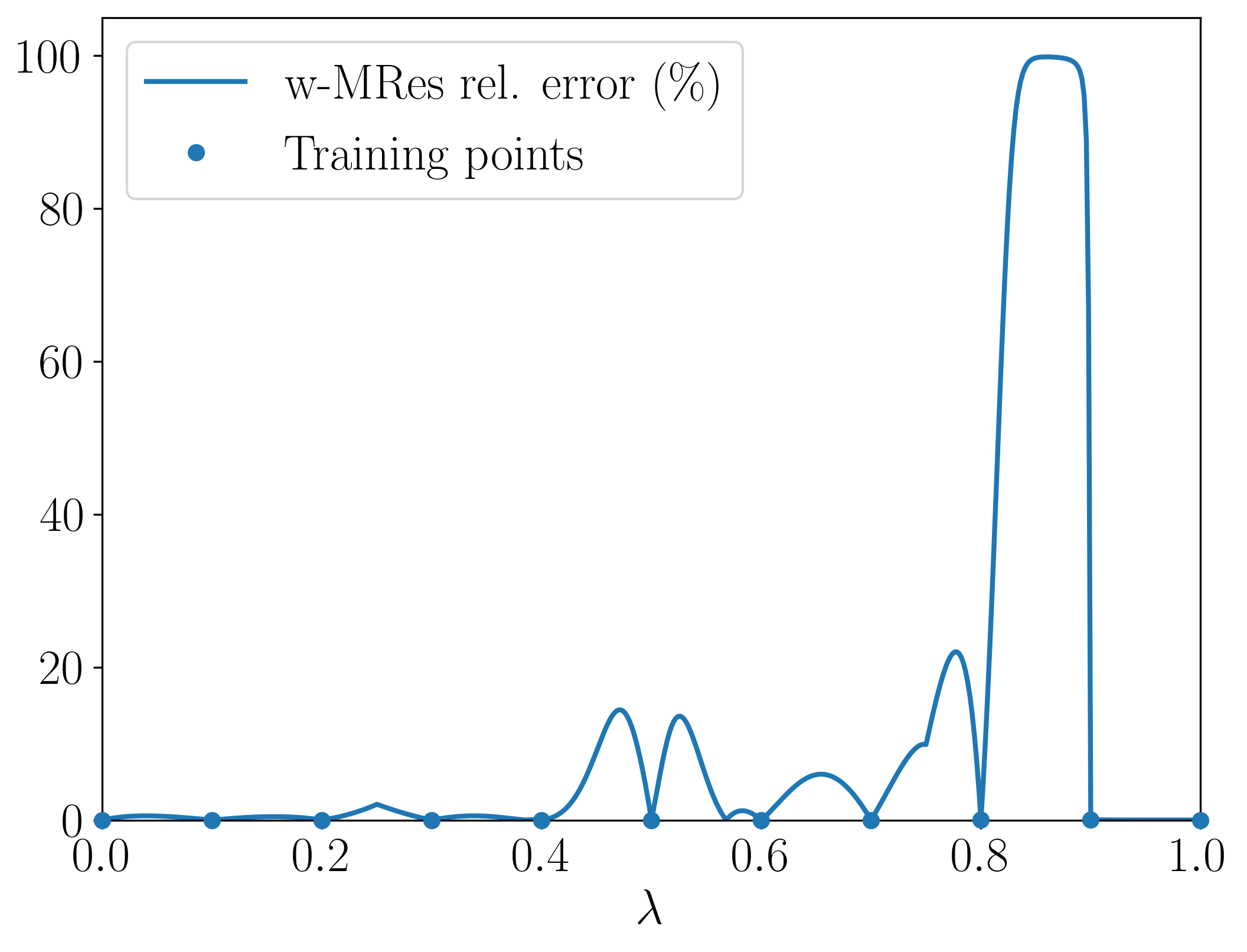}
    \caption{w-MinRes relative error.}
    \label{sfig:ex02_a}
    \end{subfigure}
    \begin{subfigure}[b]{0.34\textwidth}
    \includegraphics[width=\textwidth, height=121px]{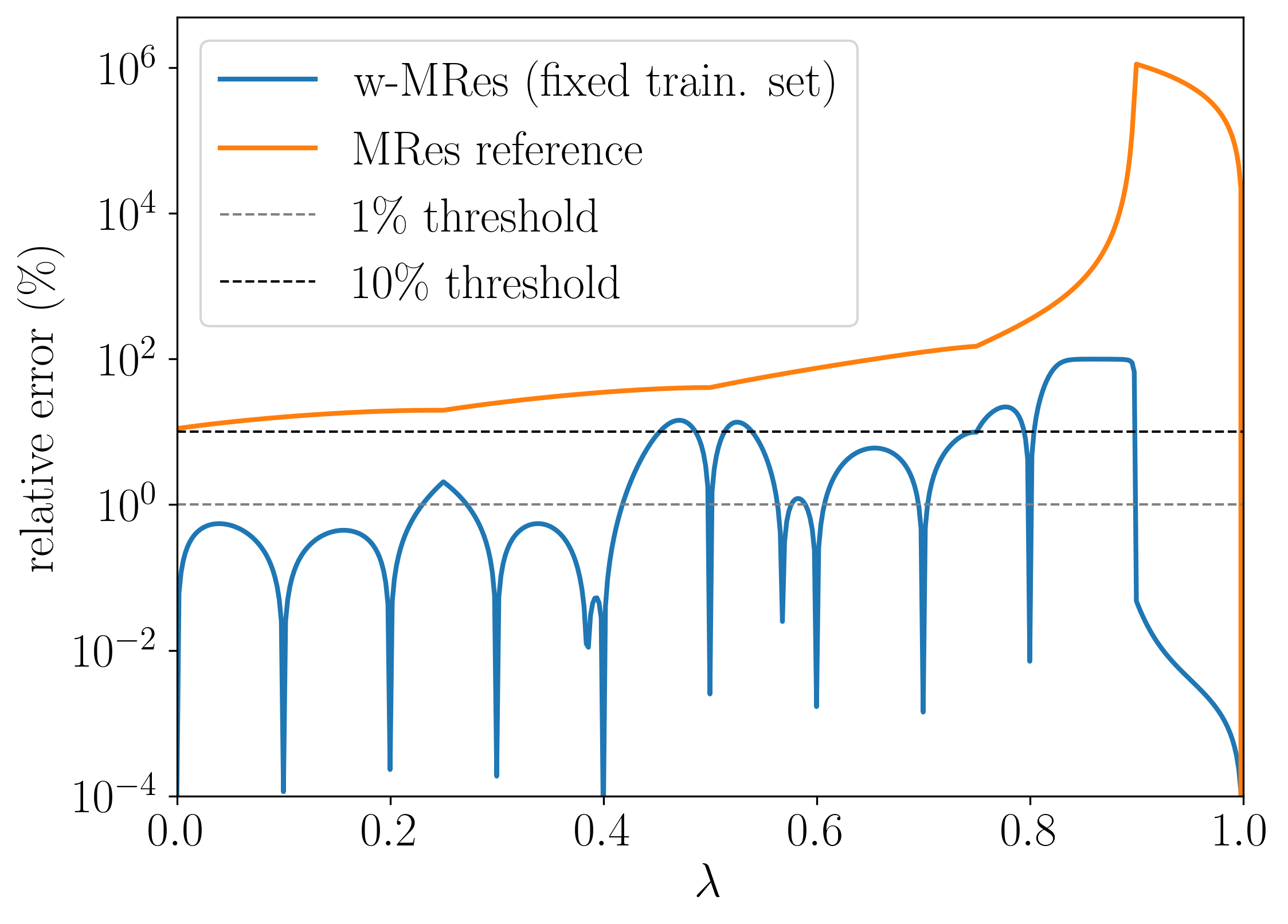}
    \caption{Relative error comparison.}
    \label{sfig:ex02_b}
    \end{subfigure}
    \begin{subfigure}[b]{0.32\textwidth}
    \includegraphics[width=\textwidth, height=121px]{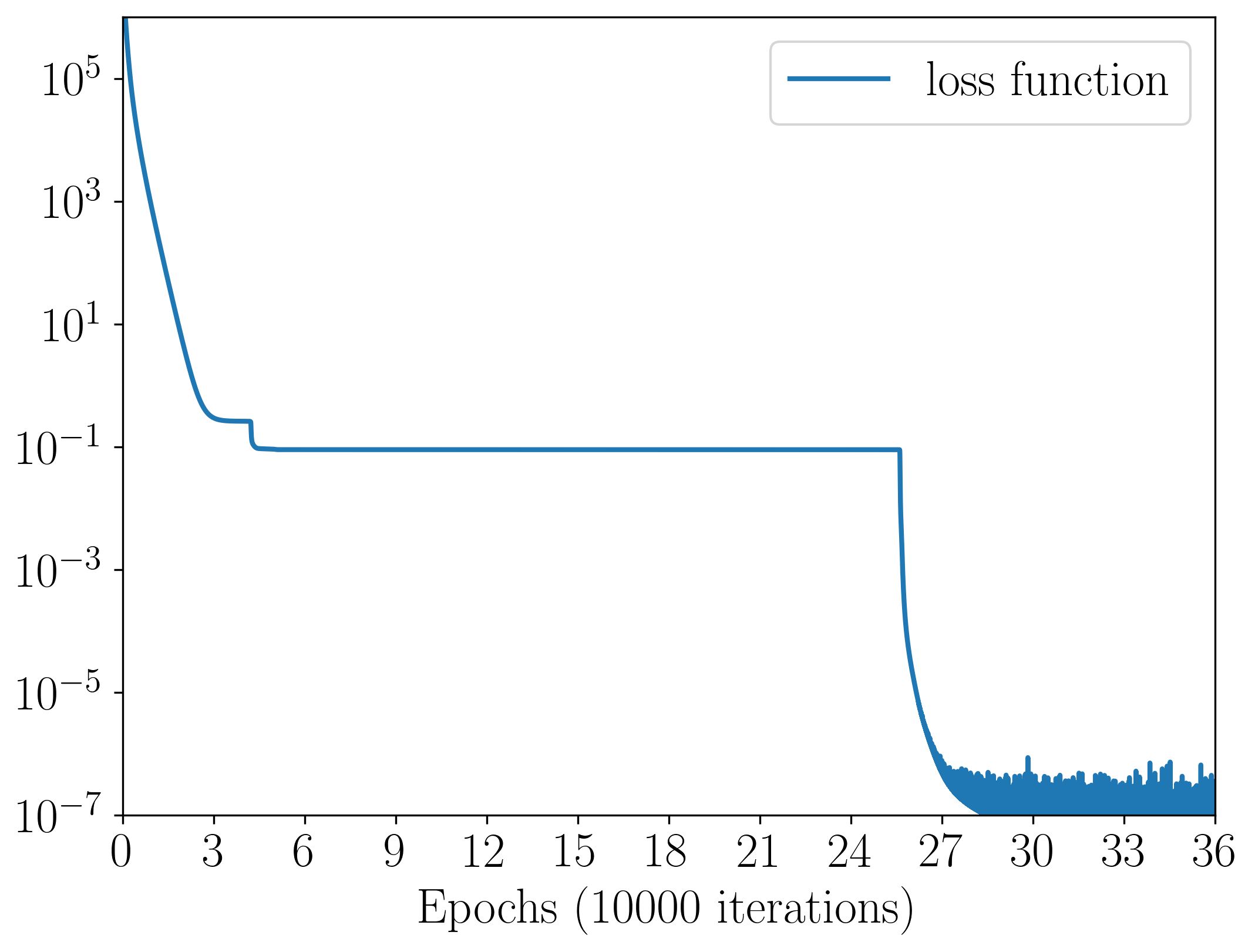}
    \caption{Loss function.}
    \label{sfig:ex02_c}
    \end{subfigure}
    \begin{subfigure}[b]{0.32\textwidth}
    \includegraphics[width=\textwidth, height=121px]{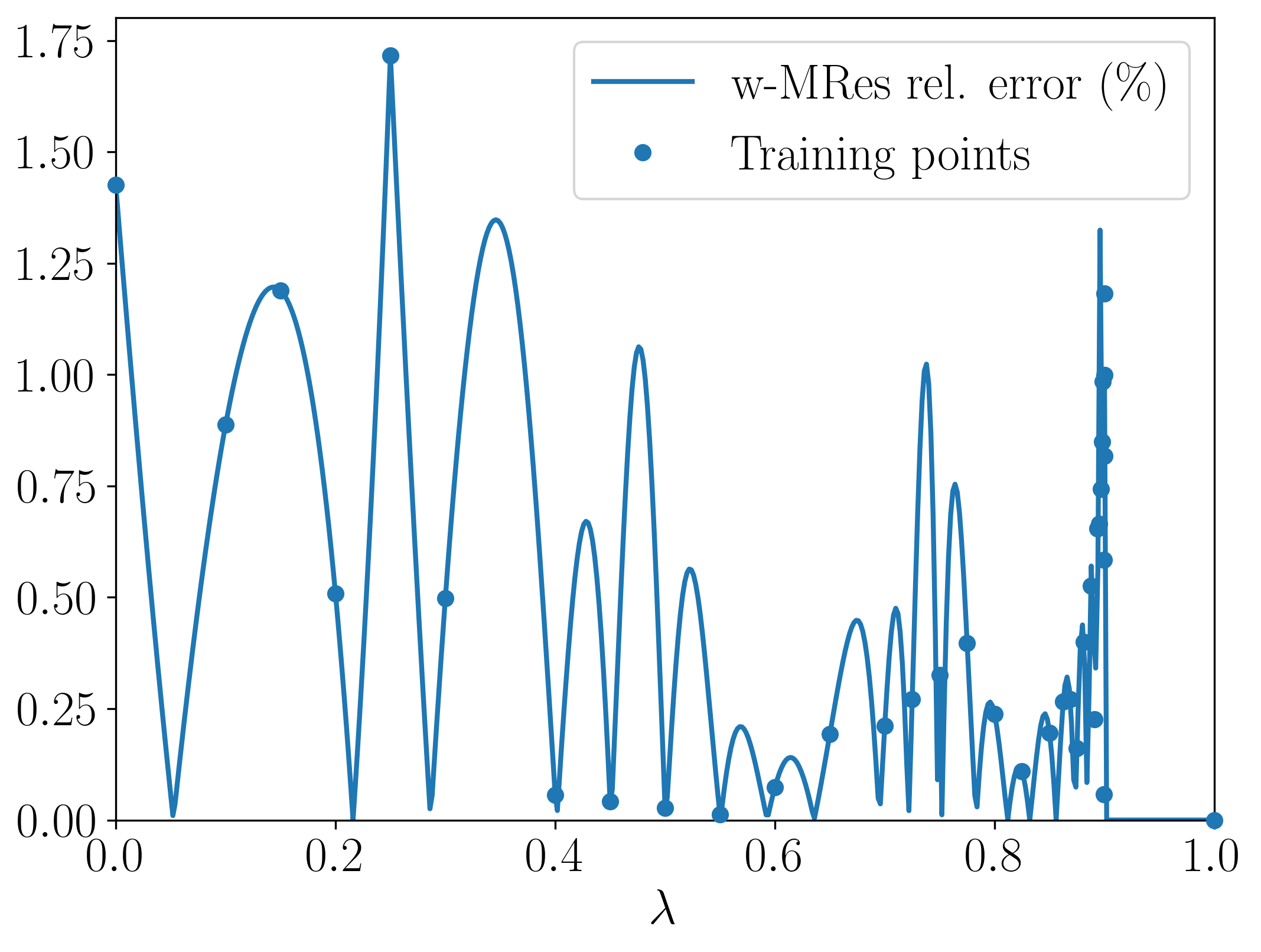}
    \caption{w-MinRes relative error.}
    \label{sfig:ex02_d}
    \end{subfigure}
    \begin{subfigure}[b]{0.34\textwidth}
    \includegraphics[width=\textwidth, height=121px]{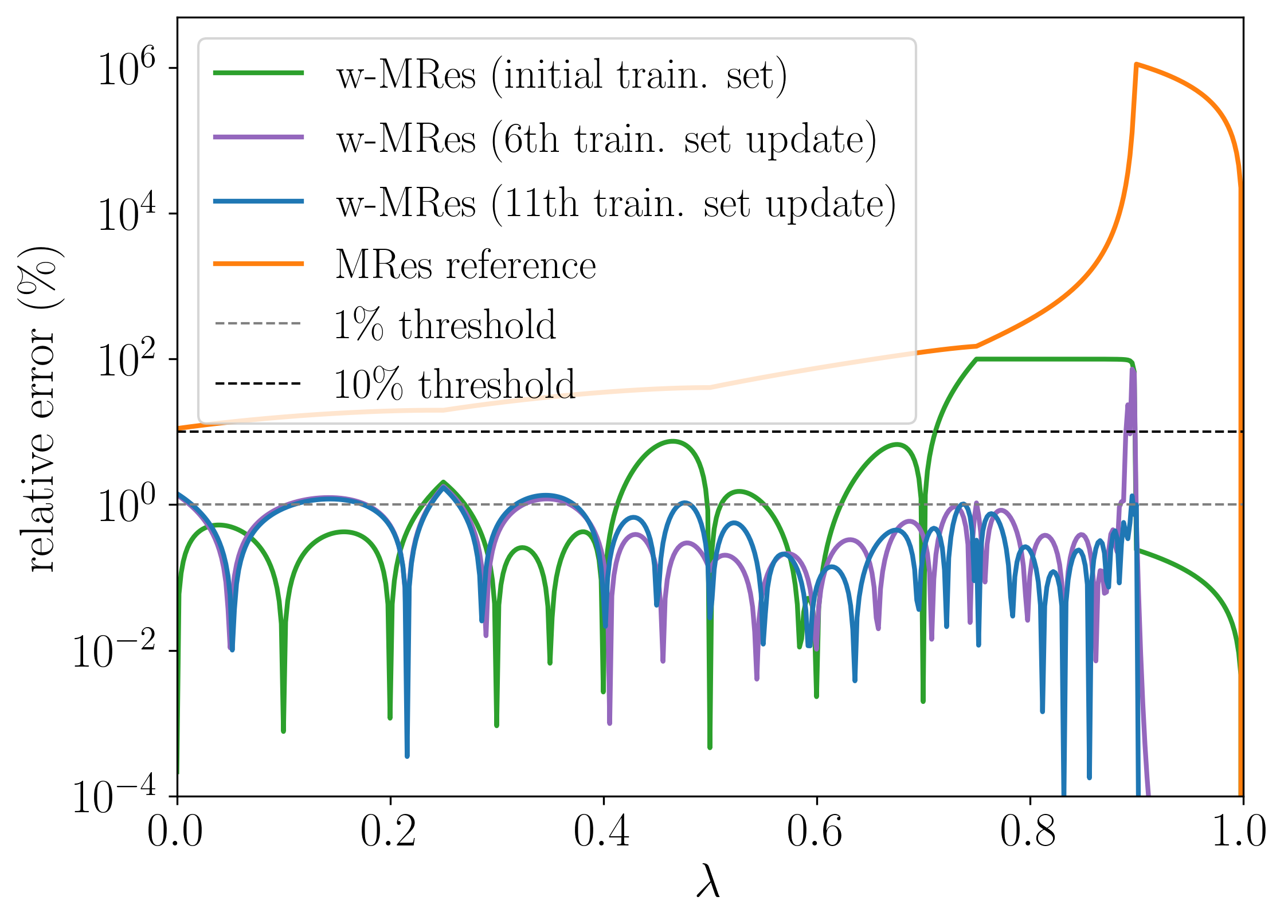}
    \caption{Relative error comparison.}
    \label{sfig:ex02_e}
    \end{subfigure}
    \begin{subfigure}[b]{0.32\textwidth}
    \includegraphics[width=\textwidth, height=121px]{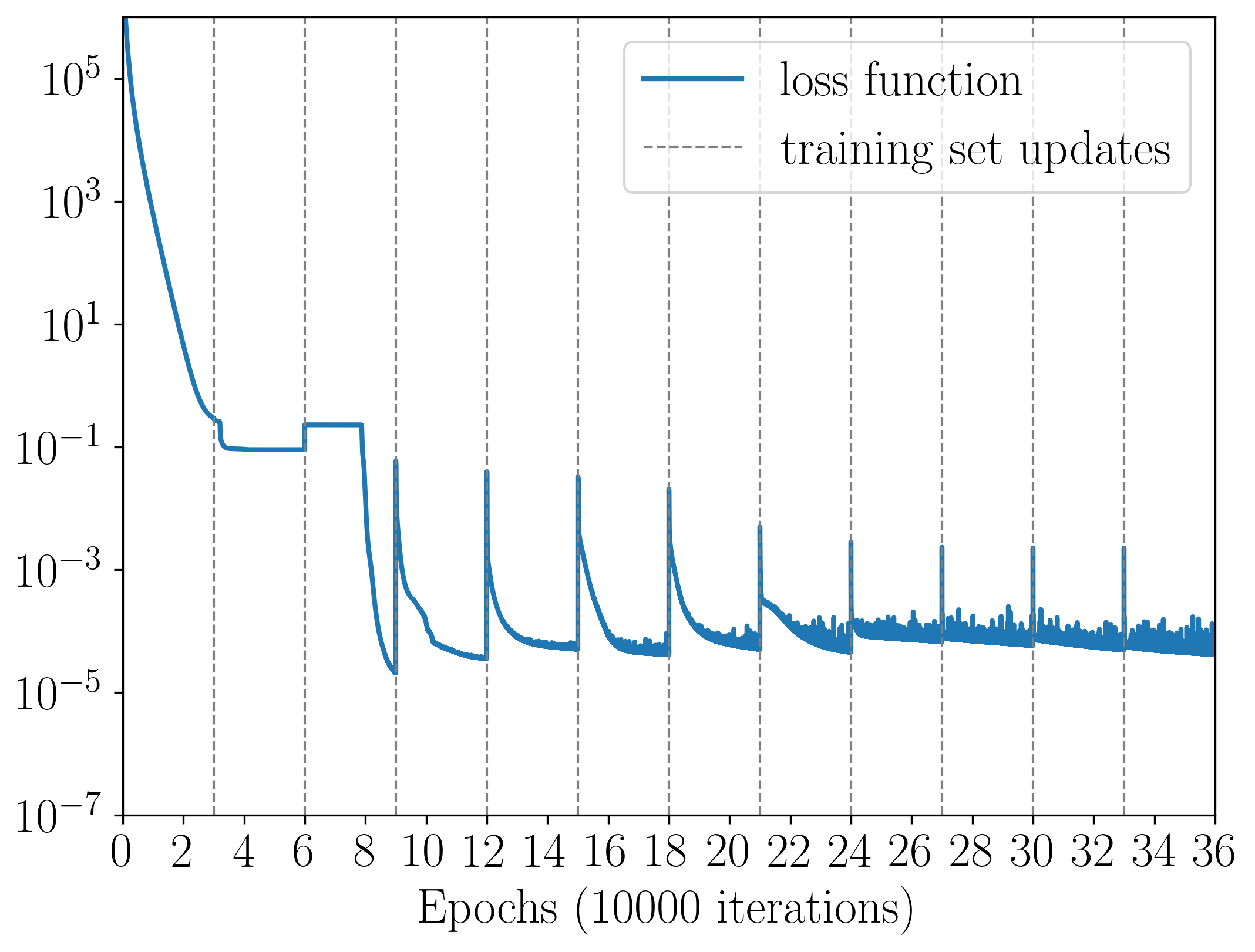}
    \caption{Loss function.}
    \label{sfig:ex02_f}
    \end{subfigure}
    \caption{ML-MinRes performance for a 1D advection example with parametric right-hand side. First row shows the results for a fixed training set. The second row displays the results for an adaptive strategy.
    }  
\label{fig:example_2_results_advection_parametric_rhs}
\end{figure}

\begin{figure}
    \centering
    \begin{subfigure}[b]{0.33\textwidth}
    \includegraphics[width=\textwidth, height=121px]{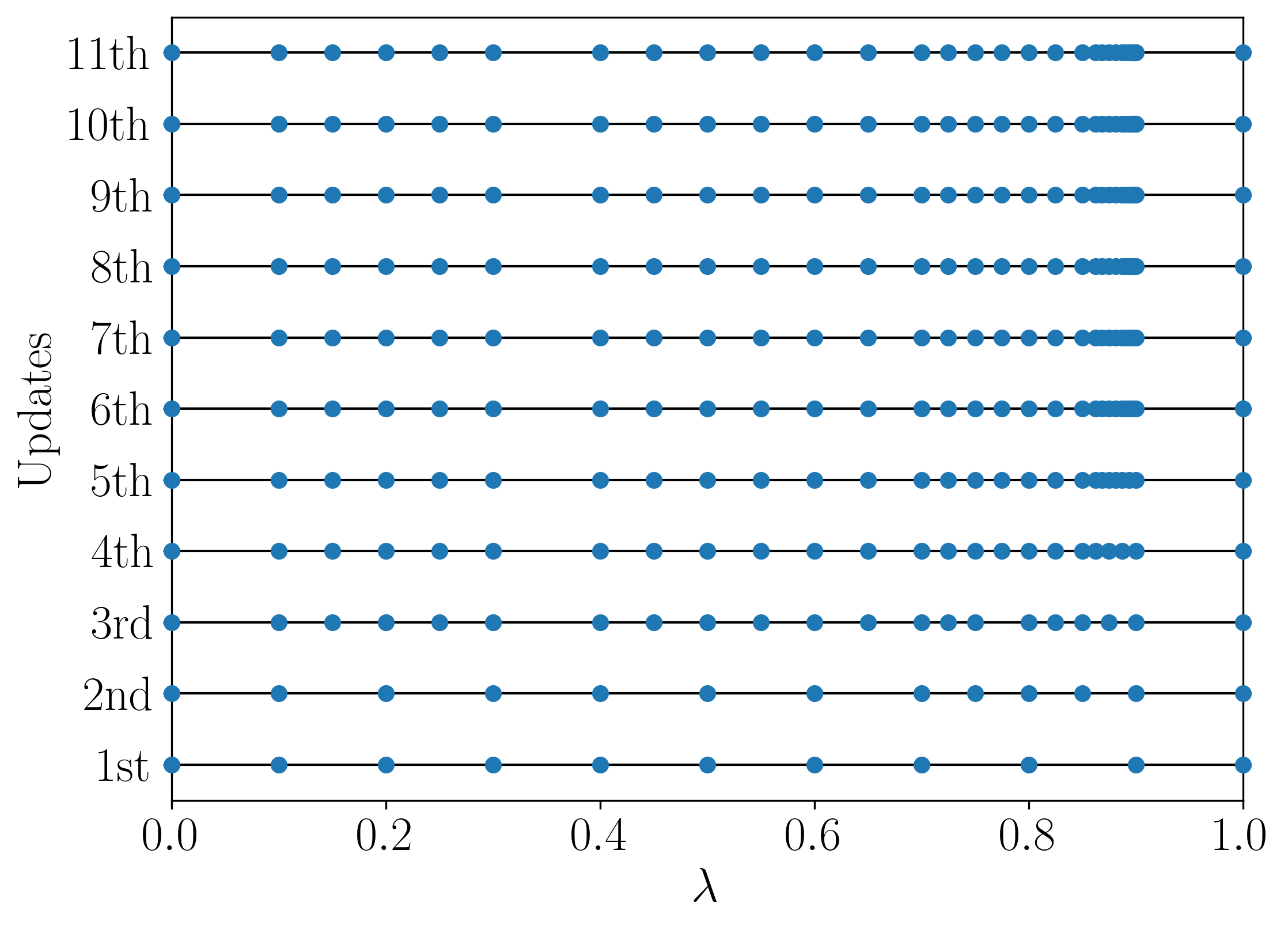}
    \caption{Adaptive training-points}
    \label{sfig:ex02_g}
    \end{subfigure}
    \begin{subfigure}[b]{0.33\textwidth}
    \includegraphics[width=\textwidth, height=121px]{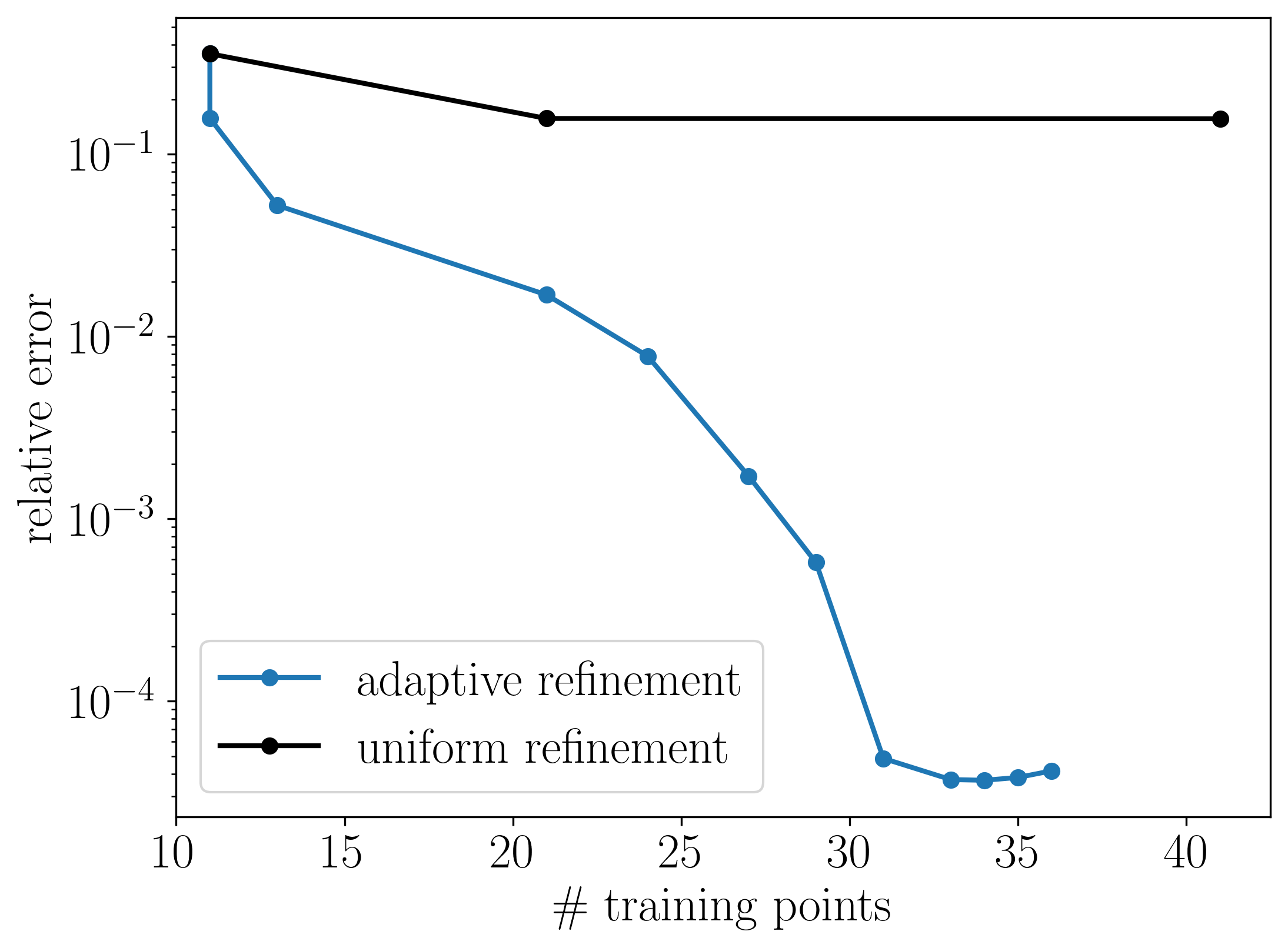}
    \caption{Refinements comparison}
    \label{sfig:ex02_h}
    \end{subfigure}    
\caption{ML-MinRes Adaptive selection of training set and a comparison with uniform refinement for a 1D advection example with parametric right-hand side.}  
\label{fig:example_2_adaptivity_advection_parametric_rhs}
\end{figure}

Figure~\ref{fig:example_2_results_advection_parametric_rhs} shows the results for the ML-MinRes with a fixed trained set (first row) and with the adaptive training set strategy (second and third rows) presented in Subsection~\ref{subsec:adaptive_training}. Subfigures in the first column show the relative error in linear scale for the training with a fixed training set (Subfigure~\ref{sfig:ex02_a}) and the adaptive training set (Subfigure~\ref{sfig:ex02_d}). Subfigures~\ref{sfig:ex02_b} and~\ref{sfig:ex02_e} compare the relative error in logarithmic scale between MinRes and ML-MinRes with fixed and adaptive training sets. The last column displays the loss function for both cases. 
Notice that the fixed training set with eleven elements poorly captures the quantity of interest behavior in this example. This is because the training set does not represent the map between $\oldlambda$ and the quantity of interest. Then, the relative error increases drastically between 0.8 and 0.9 (see Subfigure~\ref{sfig:ex02_a}), even when the loss function is small. Thus, the adaptive training set emerges as a good alternative to overcome the problem in this type of scenario.

Finally, Figure~\ref{fig:example_2_adaptivity_advection_parametric_rhs} shows the adaptive training-points generation and a comparison with uniform refinement. In Subfigure~\ref{sfig:ex02_h} we can see the behavior of the relative error for a test set (different from the adaptive training set) when we add points to the training set with adaptive (blue line) and uniform (black line) refinements after $30000$ epoch per refinement.

\subsection{2D discontinuous-media diffusion equation with two parameters}

\begin{figure}[h]
    \centering
    \begin{subfigure}[b]{0.32\textwidth}
    \includegraphics[width=\textwidth]{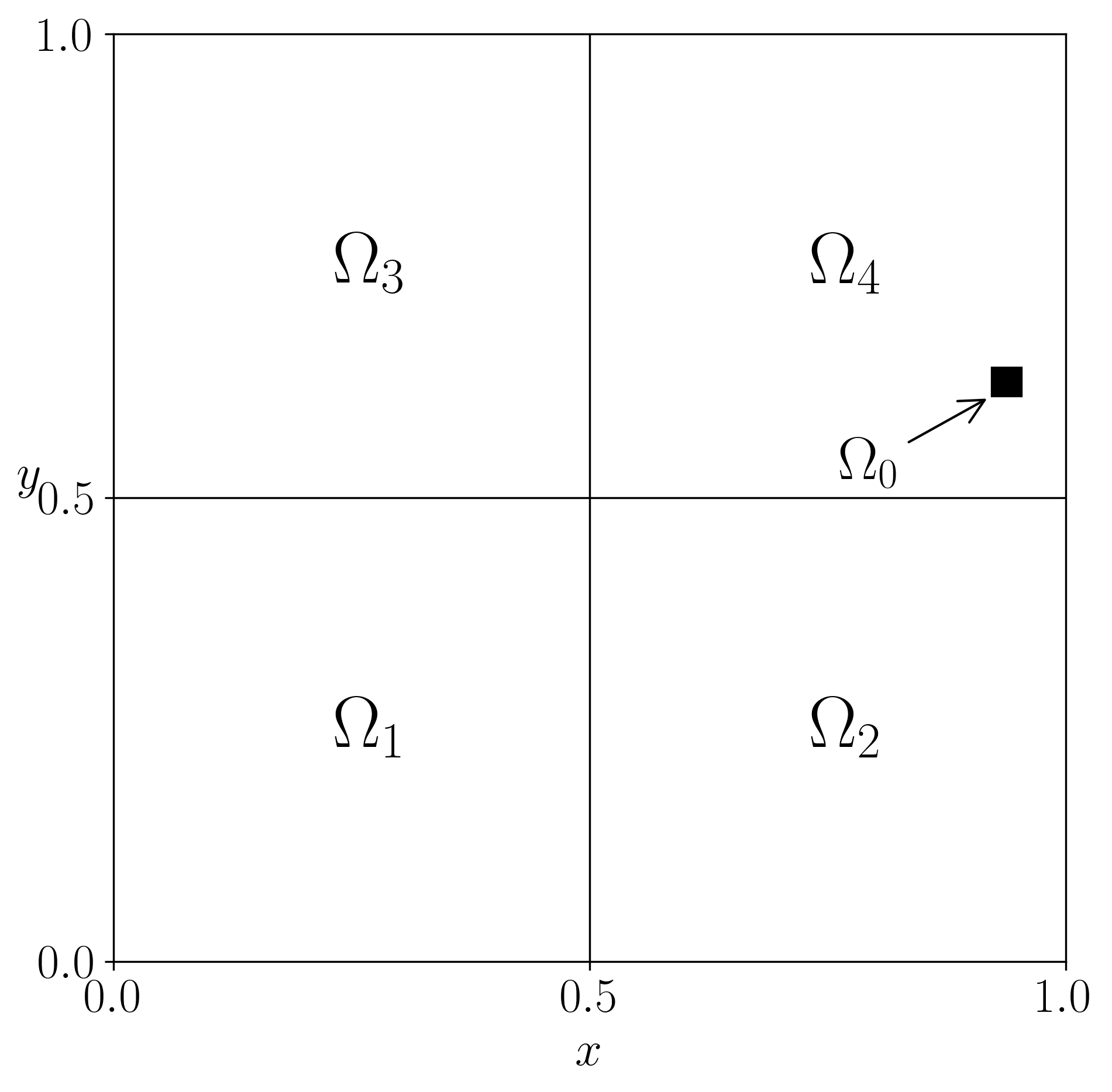}
    \caption{domain subdivision on $\Omega$}
    \label{subfig:Thermal_con_division}
    \end{subfigure}
    \begin{subfigure}[b]{0.32\textwidth}
    \includegraphics[width=\textwidth]{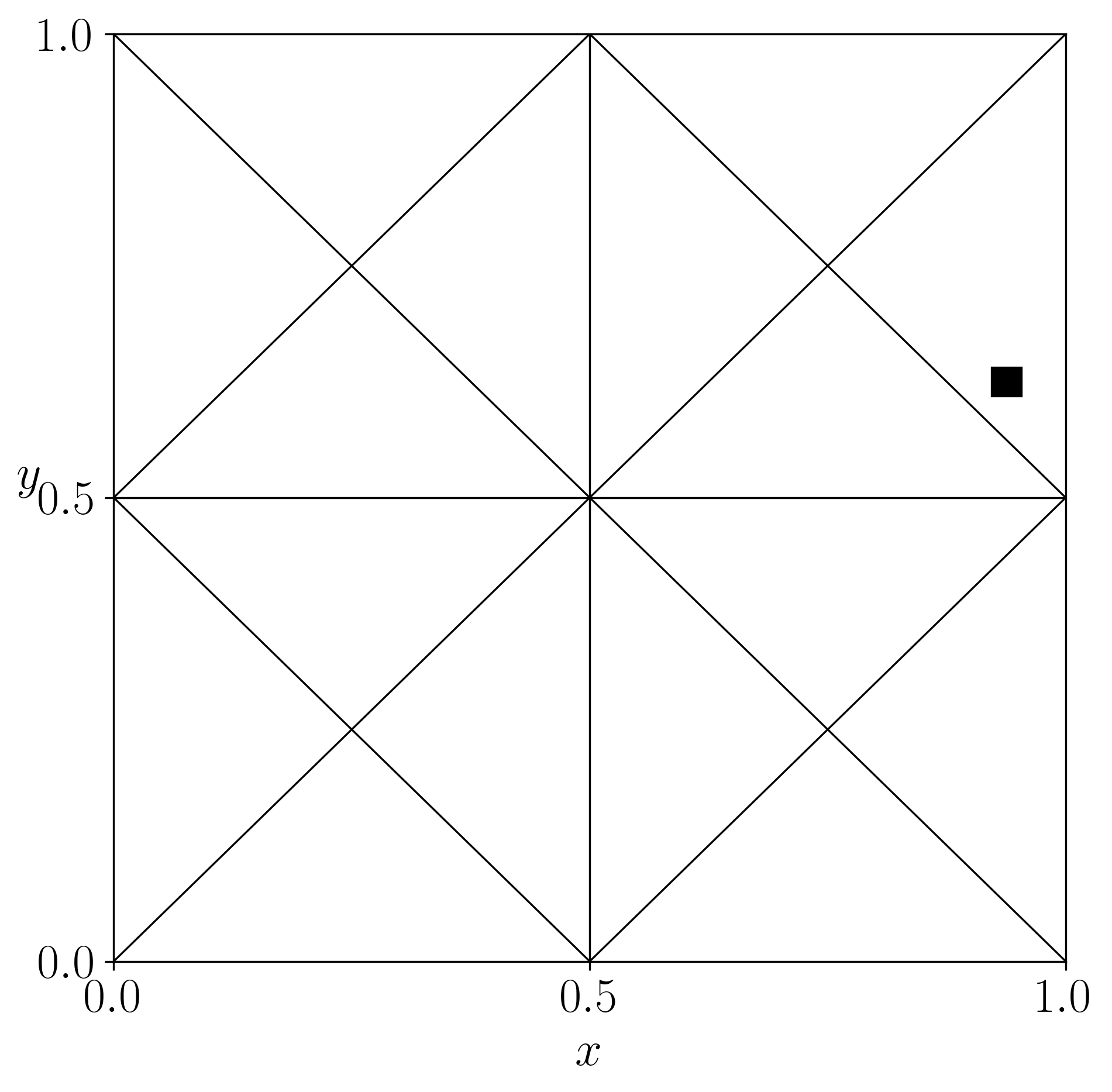}
    \caption{trial-space mesh}
    \label{fig:trial_mesh}
    \end{subfigure}
    \begin{subfigure}[b]{0.32\textwidth}
    \includegraphics[width=\textwidth]{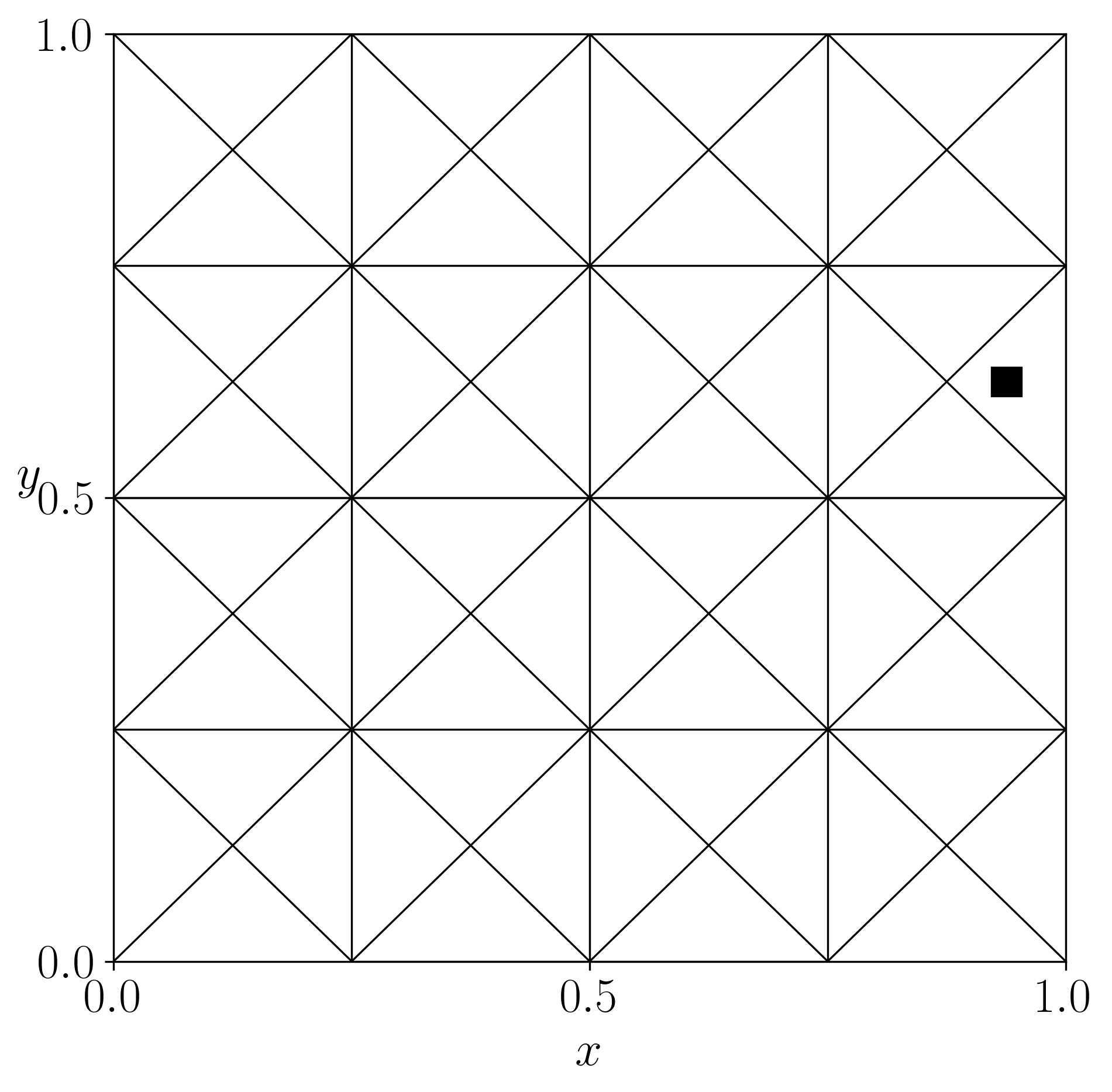}
    \caption{test-space mesh}
    \label{fig:test_mesh}
    \end{subfigure}
    \caption{Domain subdivision and meshes for discretization. Back square represents the QoI location.}
    \label{fig:ex06_domain_division}
\end{figure}

Now, we consider a 2D diffusion equation, where the diffusion coefficient is given by a two-dimensional parameter $\lambda = (\alpha,\beta)$. The partial differential equation is the following: 
\begin{equation}
\left\{
\begin{array}{rl}
- \nabla\cdot (a(x,y;\lambda) \nabla u) = 1, & \hbox{in } \Omega :=(0,1)^2, \\
u = 0, & \hbox{on } \partial \Omega,
\end{array}
\right.
\end{equation}
where
$$a(x,y;\lambda) = \begin{cases} \alpha, & \text{if } (x,y) \in \Omega_1\cup \Omega_3, \\  \beta, & \text{if } (x,y) \in \Omega_2\cup \Omega_4,
\end{cases}$$ 
with $\Omega_1 = (0 \,, 0.5)^2$, $\Omega_2 = (0.5 \,, 1)\times(0\, , 0.5)$, $\Omega_3 = (0\, , 0.5)\times (0.5\, , 1)$, and $\Omega_4 = (0.5\, , 1)^2$. See Subfigure~\ref{subfig:Thermal_con_division}.

Here, the quantity of interest will be given by 
$$q(w)= \frac{1}{|\Omega_0|}\int_{\Omega_0} w\,,\quad \forall w\in\mbbU,$$
where, $\Omega_0 = (0.9375 - \frac{1}{2^6}\, , 0.9375 + \frac{1}{2^6})\times (0.625 - \frac{1}{2^6} \, 0.625 + \frac{1}{2^6})$.

Given continuous trial and test spaces $\mbbU=\mbbV = H_{0}^{1}(\Omega)$, we use the following variational formulation:
\begin{equation}
\left\{
\begin{array}{l@{\;}lr}
\multicolumn{3}{l}{\text{Find } u(\lambda) \in \mbbU \text{ such that:}} \\
\alpha\underbrace{\displaystyle \int_{\Omega_1\cup \Omega_3} \nabla u(\lambda) \cdot \nabla v}_{\left<B_1u(\lambda),v\right>_{\mbbV^*,\mbbV}} + \beta\underbrace{ \int_{\Omega_2\cup \Omega_4} \nabla u(\lambda) \cdot \nabla v}_{\left<B_2u(\lambda),v\right>_{\mbbV^*,\mbbV}} &= \underbrace{\displaystyle\int_{\Omega} v}_{\left<\ell,v\right>_{\mbbV^*,\mbbV}} , & \forall\, v\in \mbbV.
\end{array}
\right.
\end{equation} 
Thus, we define $B(\lambda):\mathbb U\to \mathbb V^*$ as $B(\lambda):=\alpha B_1+\beta B_2$. Additionally, we define the following weighted inner product:
\begin{equation}
\nonumber 
(v,\nu)_{\omega(\lambda)} := \int_{\Omega} \omega(\lambda)\nabla v \cdot \nabla \nu\,,\quad \forall\,v,\nu \in  \mbbV,
\end{equation}
where $\omega(\lambda)$ is a \change{piecewise} constant \change{function with a constant value} on each element $T$ of the test-space mesh $\Omega_h$ (64 constants, see Figure~\ref{fig:test_mesh}).

The discrete trial space $\mbbU^n$ is given by conforming piecewise linear functions with $n = 5$ on a coarse mesh (see Figure~\ref{fig:trial_mesh}), while the discrete test space $\mbbV^m$ is given by conforming piecewise linear functions with $m = 25$ on the mesh $\Omega_h$ depicted in Figure~\ref{fig:test_mesh}.

\change{The discrete mixed system becomes:
\begin{equation}
\nonumber
\left\{ \begin{array}{l@{}l}
\text{Find } r^{m}(\lambda)\in \mbbV^{m} \hbox{ and } u^{n}(\lambda)\in\mbbU^{n} \text{ such that for all } v^{m} \in \mbbV^{m} \hbox{ and } w^{n} \in \mbbU^{n} & \\\\
\displaystyle \sum_{T\in\Omega_h}c_T(\lambda)\int_T \nabla r^{m}(\lambda) \cdot \nabla v^{m} + \alpha\int_{\Omega_1\cup \Omega_3} \nabla u^n(\lambda) \cdot \nabla v^{m} + \beta \int_{\Omega_2\cup \Omega_4} \nabla u^{n}(\lambda) \cdot \nabla v^{m} &{} = \displaystyle \int_{0}^{1} v^{m},\\\\
\displaystyle\alpha\int_{\Omega_1\cup \Omega_3} \nabla w^n \cdot \nabla r^{m}(\lambda) + \beta \int_{\Omega_2\cup \Omega_4} \nabla w^{n} \cdot \nabla r^{m}(\lambda) &{} = 0.
\end{array} \right.
\end{equation}
}

\begin{figure}[h!]
\begin{center}
  \begin{subfigure}[b]{0.32\textwidth}
    \includegraphics[width=\textwidth, height=121px]{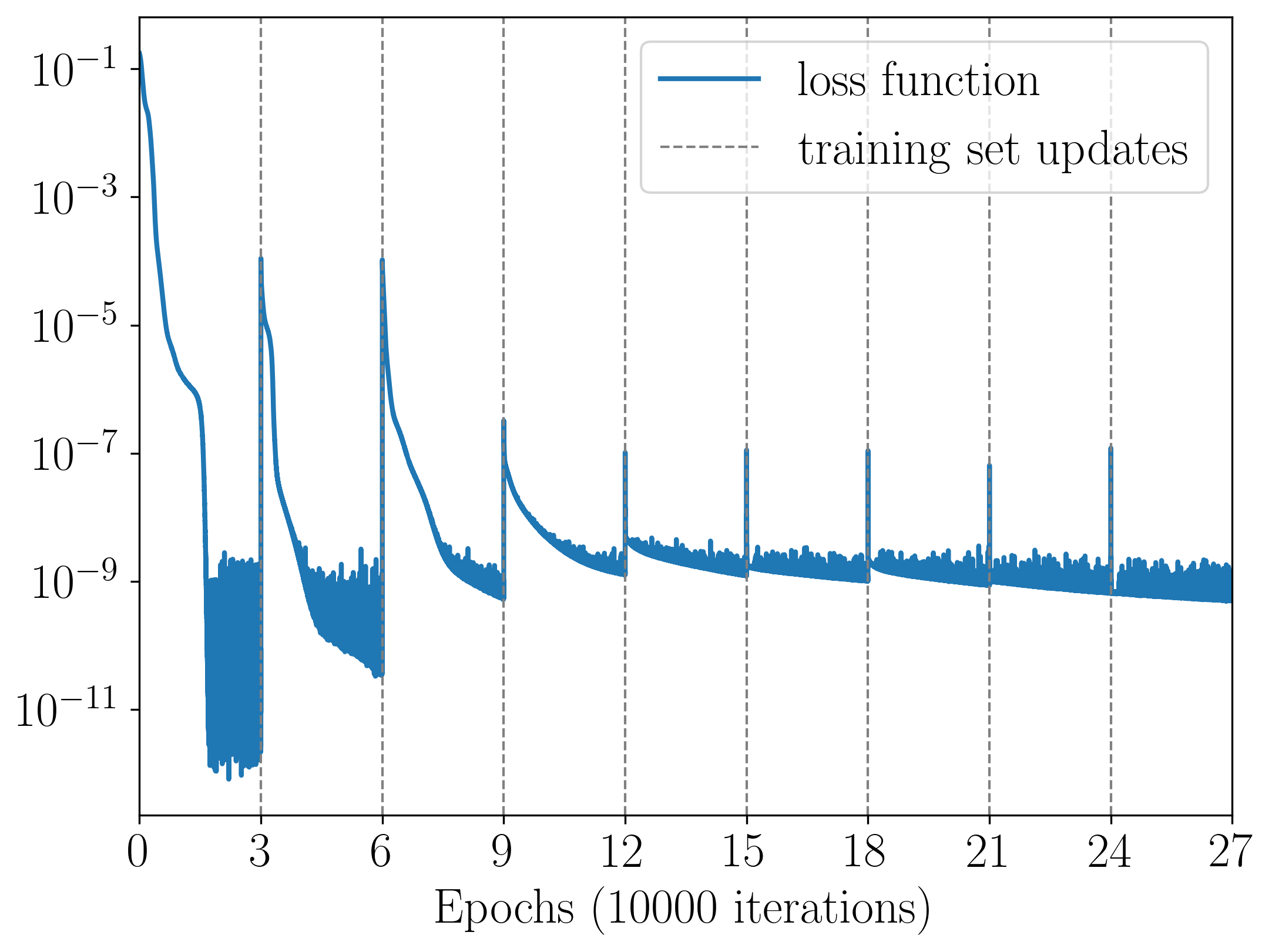}
    \caption{Loss function.}
    \label{sfig:ex07_a}
  \end{subfigure}
  \begin{subfigure}[b]{0.32\textwidth}
    \includegraphics[width=\textwidth,  height=121px]{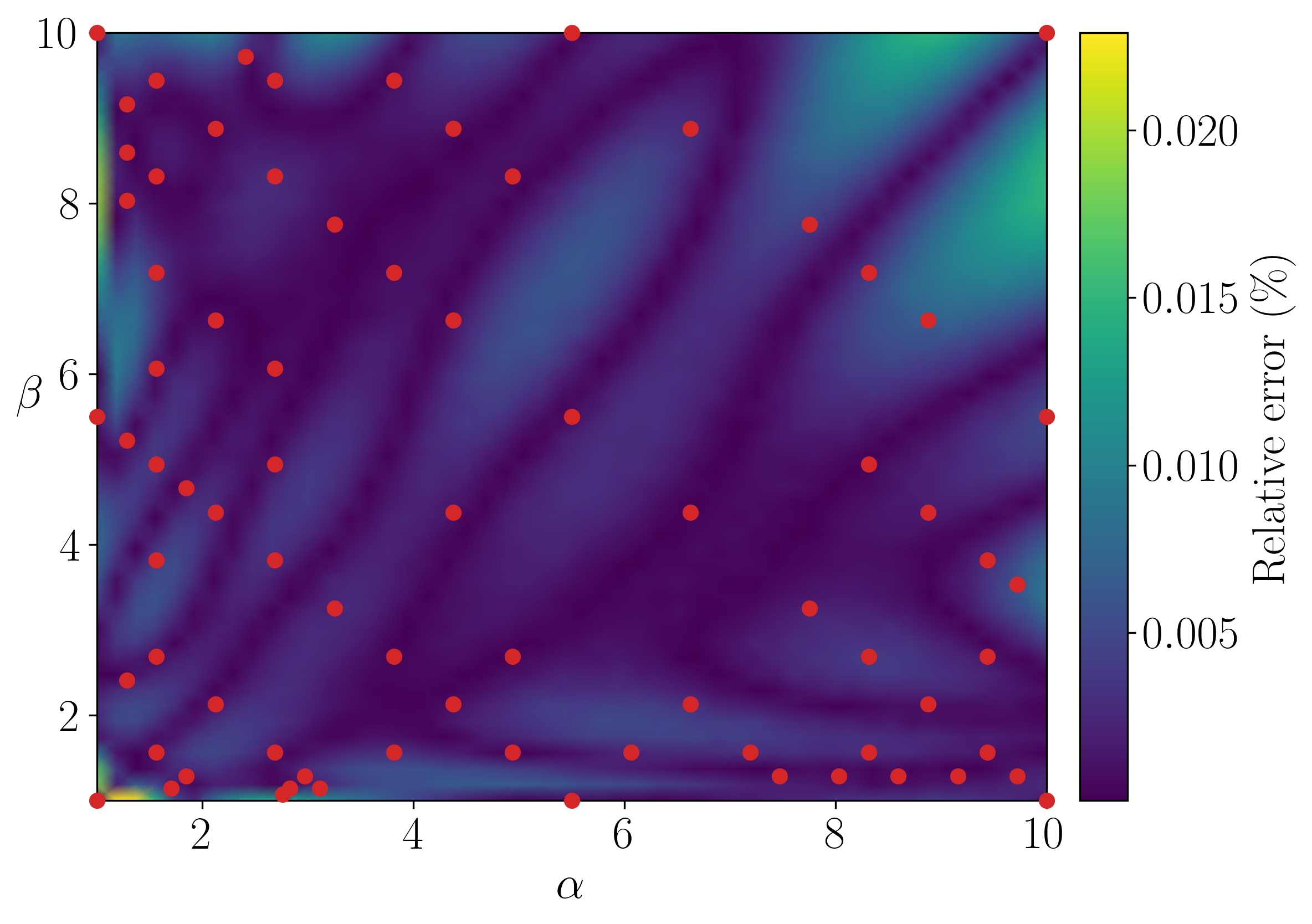}
    \caption{w-MinRes relative error.}
    \label{sfig:ex07_b}
  \end{subfigure}
  \begin{subfigure}[b]{0.32\textwidth}
    \includegraphics[width=\textwidth,  height=121px]{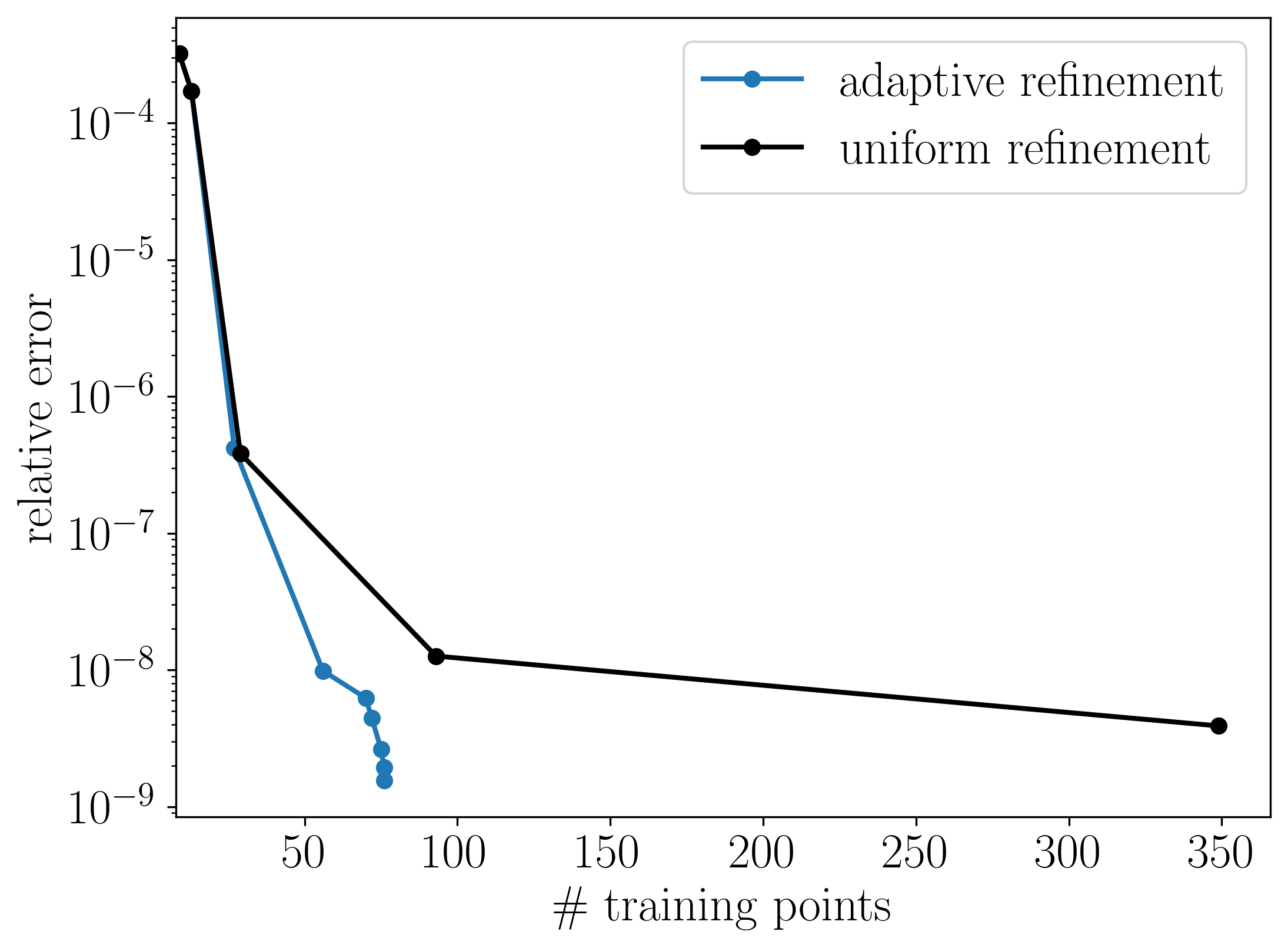}
    \caption{Refinements comparison.}
    \label{sfig:ex07_c}
  \end{subfigure}
\end{center}  
  \caption{ML-MinRes performance for two-parameters diffusion example. Red dots show the training points.}
  \label{fig:ex07_ML-MRes_results}
\end{figure}

On this example, we employ the adaptive adaptive training set with eight stages, $30000$ epoch per stage and $\gamma=5$.   
We also define the initial training set composed of $N_s = 9$ elements, composed by $\{1, 5.5, 10\}\times \{1, 5.5, 10\}$ (Cartesian product). 

Figure~\ref{fig:ex07_ML-MRes_results} displays the results for the ML-MinRes after training with an adaptive training-set strategy.
Subfigure~\ref{sfig:ex07_a} displays the loss function through the six stages of the adaptive training.
Subfigure~\ref{sfig:ex07_b} shows the the
relative error (\%) in linear scale for full test set after the training with an adaptive training set. Red dot represent the final training set.  
In Subfigure~\ref{sfig:ex07_c} we can compare the behaviour of the relative error of a test set when we add points to the training set with adaptive (blue line) and uniform (black line) refinements after $30000$ epoch per refinement.
\section{Conclusions}
\label{sec:concl}

In the present work, we extend the Machine Learning Minimal-Residual (ML-MinRes) method~\cite{BreMugZeeCAMWA2021} and \cite{NeuralControl2022} to a more general framework of parametric PDEs. The main idea is to compute some high-precision quantity of interest of the solution using finite-element coarse meshes. This finite element scheme will generate a cheap and accurate method for the quantity of interest. As part of the main new features, we highlight the following:
\begin{itemize}
    \item We avoid the problem of integrating a neural network by replacing the inner-product weight with a piecewise \change{constant} weight (with coefficients delivered by a neural network). In this way, we have exact integration with quadrature rules.
    \item The piecewise \change{constant} weight for the inner product and the affine decomposition for $B(\lambda)$ and $\ell(\lambda)$ allow to pre-assemble the matrices. This fact helps to perform more efficient training and faster online evaluation after training.  
    \item We propose an adaptive training set that allows for faster and more accurate neural network training. From the numerical experiments, we can see that adaptive training performs better than uniform refinements. In addition, it helps to reduce overfitting.  
\end{itemize}

\section*{Acknowledgements}
\label{sec:acknow}

\change{The authors are grateful to the reviewer whose comments have improved the accessibility of our work to a larger audience.}
This publication has received funding from the European Union's Horizon 2020 research and innovation programme under the Marie Sklodowska-Curie grant agreement No 777778 (MATHROCKS). 
IB has also received funding from ANID FONDECYT/Postdoctorado No 3200827 and the Engineering and Physical Sciences Research Council (EPSRC), UK under Grant EP/W010011/1. 
The work by IM was done in the framework of the Chilean grant ANID FONDECYT \#1230091.
DP has received funding from: the Spanish Ministry of Science and Innovation projects with references TED2021-132783B-I00, PID2019-108111RB-I00 (FEDER/AEI) and PDC2021-121093-I00 (MCIN / AEI / 10.13039/501100011033/Next Generation EU), the ``BCAM Severo Ochoa'' accreditation of excellence CEX2021-001142-S / MICIN / AEI / \\ 10.13039/501100011033; the Spanish Ministry of Economic and Digital Transformation with Misiones Project IA4TES (MIA.2021.M04.008 / NextGenerationEU PRTR); and the Basque Government through the BERC 2022-2025 program, the Elkartek project SIGZE (KK-2021/00095), and the Consolidated Research Group MATHMODE (IT1456-22) given by the Department of Education. 
KvdZ was supported by the Engineering and Physical Sciences Research Council (EPSRC), UK under Grant EP/T005157/1 and EP/W010011/1.
\newpage
\appendix
\bibliography{bibliography}
\bibliographystyle{siam}
\end{document}